\documentclass[fleqn,12pt,twoside]{article}
\usepackage{espcrc1}
\usepackage{a4wide}
\usepackage{graphicx,palatino,pifont,times}
\usepackage{oldgerm}
\usepackage{amsmath,amsthm}
\usepackage{amssymb}
\usepackage{amstext}
\usepackage{epsfig}
\usepackage{psfig}
\usepackage{cite}

\newtheorem{thm}{Theorem}[section]
\newtheorem{Def}{Definition}[section]
\newtheorem{Rem}{Remark}[section]
\newtheorem{prop}{Proposition}[section]
\newtheorem{corl}{Corollary}[section]

\newtheorem{lm}{Lemma}[section]

\def\N{{I\!\!N}}
\def\CC{{\rm \kern.24em \vrule width.02em height1.4ex depth-.05ex \kern-.26emC}}
\def\TagOnRight
\def\AA{{it I} \hskip-3pt{\tt A}}

\def\QQ{\rlap {\raise 0.4ex \hbox{$\scriptscriptstyle |$}} {\hskip -0.1em Q}}

\catcode`\@=11
\def\theequation{\@arabic{\c@section}.\@arabic{\c@equation}}
\catcode`\@=12
\begin{document}
\begin{center}
{\Large \bf { Coalescence Hidden Variable Fractal Interpolation
Functions and its Smoothness Analysis}}
\end{center}
\begin{center}
{\sc{A.K.B.Chand\footnote{ The present research is partially supported by 
CSIR Grant No: 9/92(160)/98-EMR-I, India.}and G.P.Kapoor}}\\
Department of Mathematics,\\
Indian Institute of Technology Kanpur,\\
Kanpur 208016,India \\
Email: chand@iitk.ac.in ; gp@iitk.ac.in\\
Phone: 91-512-2597609\\
Fax : 91-512-2597500\\
\end{center} 
\begin{abstract}
{\bf ABSTRACT :} We construct a coalescence  hidden variable  fractal
interpolation function(CHFIF) through a non-diagonal iterated function
system(IFS). Such a FIF may be self-affine or non-self-affine depending
on the parameters of the defining non-diagonal IFS. The smoothness analysis of the CHFIF 
has been carried out by using the operator approximation technique. The deterministic
construction of functions having order of modulus continuity 
$O (|t|^{\delta} (\log|t|)^m)$ ($m$ a non-negative integer and
$ 0 < \delta \le 1$)  is possible through our CHFIF. The bounds of fractal dimension
of  CHFIFs are obtained first  in certain  critical cases
and then, using estimation of these bounds, 
the bounds  of fractal dimension of any FIF are found.
\end{abstract}
\section{INTRODUCTION}
The fractal curves arise during several applications in  various
disciplines such as Natural Science\cite{M,IK,MTA,SGV},
Engineering Applications\cite{VLT}, Economics\cite{K} etc. 
To approximate these curves,  Barnsley\cite{B,B1} constructed
a fractal interpolating function (FIF) arising from a suitable
iterated function system (IFS).
FIFs are generally self-affine in nature and the Hausdorff-Besicovitch
dimensions of their graphs are non-integers. To approximate non-self-affine
patterns, the hidden variable FIFs (HFIFs)  are constructed in \cite{B1,B2,PR} by projection
of vector valued FIF from generalized interpolation data.
However, in practical applications of FIF, the interpolation data might be
generated simultaneously from self-affine and non-self-affine functions.
Thus, the question  whether it is possible to construct an IFS that is capable of
generating both of the self-affine or non-self-affine FIFs simultaneously needs to be settled.
The hidden variable bivariate fractal interpolation surfaces
are studied in \cite{CK} by introducing the concept of {\it constrained free variables}.
In the present work a {\it Coalescence Hidden Variable FIF} (CHFIF) that
is self-affine or non-self-affine depending on the parameters of
defining IFS is constructed.
\par
Since FIFs are continuous but generally nowhere differentiable functions,
their analysis can not be done satisfactorily by restricting to classical analytic tools.
For the applications of FIF theory, in general, an expansion of the FIF
in terms of a suitable function system is usually considered.
 Barnsley and Harrington~\cite{BH} used  shifted
composition to express  affine FIFs and computed their fractal
dimensions. However, this representation is  somewhat difficult to use.
Zhen\cite{ZS} gave  another series representation of self-affine FIF through
a new function $\psi_{\sigma} \kappa_{\omega}$ to study the
H\"{o}lder property of FIF. Since, the function $\psi_{\sigma} \kappa_{\omega}$
has too many points of discontinuity, it is slightly tedious to analyze it in applications.
Zhen and Gang\cite{SG} expanded equidistant FIF on $[0,1]$ by using
Haar-wavelet function system and obtained  their global H\"{o}lder property,
when  the number of interpolation points is $N = 2^p + 1$, $p$ being a definite 
positive integer. Gang\cite{G} employed the technique of
operator approximation to characterize the H\"{o}lder
continuity of self-affine FIFs on a general set of nodes on $[0,1]$.
Bedford\cite{TB} obtained the H\"{o}lder exponent $h$ of a self-affine
fractal function that has non-linear scaling, using code space
of $n$ symbols associated with the IFS. He also showed the existence
of a larger H\"{o}lder exponent $h_{\lambda}$ defined at almost
every point with respect to Lebesgue measure. The distribution of points where the FIF has
strongest singularity is found by Maslyuk~\cite{MASLYUK}
that helps in calculating the parameters of an IFS with aid of
wavelet-based techniques, such as modulus maxima lines tracing.
The H\"{o}lder exponent needed in smoothness analysis of non-self-affine
FIF is not yet studied due to interdependence of the components
of vector valued FIF in the construction of HFIFs. 
\par
It is seen in the present paper that, contrary to the observation  of  
Barnsley\cite{B1} that `the graph of HFIF is not self-similar or self-affine or self-anything',
CHFIF is indeed self-affine under certain conditions even though the class of CHFIFs is a
subclass of the class of HFIFs. Our approximation of CHFIF  is obtained 
through an operator found with integral averages on each subinterval of the FIF. Using this
approximation,  the   H\"{o}lder exponent of the non-self-affine functions 
 arising from IFS is found for the first time.
The bounds of Fractal dimension of the CHFIF in
critical cases obtained in the present paper help to calculate the bounds
of Fractal dimension of any FIF by converting the CHFIF to a self-affine FIF.
\par
The organization of the paper is as follows: In Section 2,  
we construct a coalescence  hidden variable FIF. For this
 purpose, an IFS is constructed in $\mathbb{R}^3$ with the introduction
of constrained free variable. The projection of the attractor
of our IFS on $\mathbb{R}^2$ is a  CHFIF or a self-affine FIF depending upon  choices of hidden variables. 
The H\"{o}lder continuity of CHFIFs (both self-affine and non-self-affine) is investigated in 
Section 3 by using the operator approximation technique.
The bounds on fractal dimension of CHFIFs in
critical cases are obtained in Section 4. The results found in
the present work through Sections 2-4 are illustrated in Section 5 with 
the help of suitably chosen examples.
\section{CONSTRUCTION OF  CHFIF}
\setcounter{equation}{0}
\subsection{ Construction of IFS for CHFIF }
Let the interpolation data be  $\{(x_i,y_i) \in \mathbb{R}^2 : i= 0,1,2,\dots,N \},$
where $-\infty < x_0 < x_1 < \dots <x_N < \infty$. 
For constructing an interpolation function
$ f_1: [x_0, x_N] \rightarrow \mathbb{R} $ such that $f_1(x_i) = y_i \;\text{for all} \;
i=0,1,2,\dots N,$  consider a generalized
set of data  $\{(x_i,y_i,z_i) \in \mathbb{R}^3 | i=0,1,2,\dots,N \} $,
where $ z_{i},\; i=0,1,2,\dots,N $ are  real parameters.
The following notations are used throughout the sequel:
$I = [x_{0},x_{N}], \; I_i = [x_{i-1},x_i], \;
  g_{1} =\underset{i}{Min}  \;y_i, \; g_{2} =
\underset{i}{Max} \; y_{i}, \;
h_{1} =\underset{i}{Min} \; z_{i}, \;  h_{2} =
\underset{i}{Max} \; z_{i}$ and
$K= I \times D,$ where $D = J_1 \times J_2,$ $J_1, J_2$ are suitable compact
sets in $\mathbb{R}$ such that $[g_1, g_2] \times [h_1, h_2] \subset D.$ 
Let $L_i : I \longrightarrow I_i$ be a contractive homeomorphism and  $ F_i : K
\longrightarrow D $ be a continuous vector valued function such that
\begin{equation}\left.\label{chf1}
\begin{split}
L_i(x_0) = x_{i-1}, &L_i(x_N) = x_i \\
F_i(x_0,y_0,z_0)   = & (y_{i-1} , z_{i-1}), \;
F_i(x_N,y_N,z_N)   =  (y_i , z_i)
\end{split}\right\} 
\end{equation}
and 
\begin{equation}\left.
\begin{split}
d(F_i(x,y,z), F_i(x^*,y,z))
&\le  c \;  |x -x^*| \\
 d(F_i(x,y,z), F_i(x,y^*,z^*))
& \le  s \; d_E((y,z), (y^*,z^*))
\end{split}\right\}
\end{equation}\label{chf2}
for all  $ i = 1,2,\dots,N$ where, $c$ and $s$ are  positive constants  with $0\le\ s <1,$  
$ (x,y,z), \; (x^*,y,z),$ $ (x,y^*,z^*) \in K $, $d$ is the sup. metric on $K$ and
$d_E$ is the Euclidean metric on $\mathbb{R}^2.$ For defining the required CHFIF,
the functions $L_i$ and  $F_i$ are chosen to be of the form $L_i(x) = a_ix + b_i $ and 
\begin{equation}\label{chf3}
F_{i}(x,y,z) = A_i (y,z)^T + (p_i(x), q_i(x))^T
\end{equation}
where, $A_i$ is an upper triangular matrix  
$\begin{pmatrix} \alpha_i & \beta_i \\
0 & \gamma_i \end{pmatrix}$ and $p_i(x)$, $q_i(x)$ are  continuous functions having at least
two unknowns. We  choose $\alpha_i$ as free variable with
$ |\alpha_i| < 1 $ and $\beta_i$ as {\it{constrained  free variable}} with respect to
$\gamma_i$ such that $ |\beta_i| +| \gamma_i| < 1$.
The generalized IFS that is needed for construction of CHFIF corresponding to the data
$\{(x_i,y_i,z_i) | $ $ i=0,1,\dots,N \}$ is now defined as
\begin{equation}\label{chf4}
\{\mathbb{R}^3 ; \omega_{i}(x,y,z)= (L_{i}(x), F_{i}(x,y,z)),\; i=1,2,\dots,N \}.
\end{equation}
It is shown in the sequel that projection of the attractor of IFS (\ref{chf4}) on $\mathbb{R}^2$
 is the desired CHFIF.
\subsection{ Existence and Uniqueness of CHFIF}
It is known\cite{B2} that the IFS defined in (\ref{chf4}) associated with the data $ \{(x_i,y_i,z_i),
 \; i=0,1,\dots,N\}$ is hyperbolic with respect to a metric $d^*$ on $\mathbb{R}^3$
equivalent to the Euclidean metric. In particular, there exists a unique nonempty
compact set $G\subseteq \mathbb{R}^3$ such that
\begin{equation}\label{attractor}
G = \underset{i=1}{\overset{N}\bigcup} \omega_{i}(G)
\end{equation}
The following proposition gives the existence of a unique vector valued function $f$ that interpolates
the generalized interpolation data and also establishes that the graph of $f$ equals  the attractor $G$
 of the generalized IFS:
\begin{prop}\label{prop2} 
The attractor $G$ (c.f. (\ref{attractor})) of the IFS  defined in (\ref{chf4})
is  the  graph of the continuous  vector valued  function  
$ f : I \longrightarrow D $  such  that  $ f (x_i)  = 
 ( y_i,z_i)$ for all $ i = 1,2,\dots,N$ i.e. $G = \{(x,y,z) : x \in I\; \text{and} \; 
f(x) = (y(x), z(x))\}.$
\end{prop}
\begin{proof} Consider the family of functions,
 $\mathcal{F} = \{ f :I \longrightarrow D |f\text{ is continuous},
 \; f(x_0) = (y_0,z_0),$ $ f(x_N) = (y_N,z_N) \}. $
For $f \;\text{and} \; g  \in \mathcal{F},$ define the metric  $\rho ( f, g ) =
\underset{x \in I} \sup \| f(x) - g(x)\|$
where,  $\| . \|$ denotes the Euclidean norm on $\mathbb{R}^2$.
 Then, $(\mathcal{F} , \rho )$ is a complete metric space.
Now, for  $x \in I_{i},$ define Read-Bajraktarevi\'{c}  operator $T$
on $(\mathcal{F} , \rho )$ as
\begin{equation}\label{chf6}
(Tf) (x) = F_i (L_i^{-1}(x), y(L_i^{-1}(x)), z(L_i^{-1}(x)))
\end{equation}
For $f \in  \mathcal{F},$ using (\ref{chf1}),
$(Tf) (x_0) = F_{1}(L_1^{-1}(x_0) , y(L_1^{-1}(x_0)) , z(L_1^{-1}(x_0)))
=F_{1} (x_0, y_0, z_0) = (y_0 , z_0).$
Similarly, $(Tf) (x_N) = (y_N , z_N).$ The function
 $Tf$ is clearly  continuous on each of the subinterval $(x_{i-1},x_i)$ for
$i = 1,2,\dots,N.$ Also, from (\ref{chf1}), it follows that
$Tf(x_i^-) = Tf(x_i^+)$ for each $i$. Consequently, $TF$ is continuous on $I.$ Thus,
$ Tf \in \mathcal{F}.$ This proves that $T$ maps $\mathcal{F}$ into itself .
\par
Next, we prove that  $T$ is a contraction map on $\mathcal{F}.$ For $f \in \mathcal{F},$
 define $y_f(x), z_f(x)$ as the $y$-value and $z$-value of the vector valued function
$Tf$ at $x$.
Let $ f, g \in \mathcal{F}\; \text{and} \;x \in I_i. $ Then,
\begin{equation*}
\begin{split}
 \rho  (Tf,Tg) = & \underset{x \in I} \sup \{ \| Tf(x) - Tg(x) \| \} \\
 = &  \underset{x \in I_i} \sup \{ \| \alpha_i (y_f(L_i^{-1}(x))
 - y_g(L_i^{-1}(x))) + \beta_i (z_f(L_i^{-1}(x)) - z_g(L_i^{-1}(x))),\\
& \gamma_i  (z_f(L_i^{-1}(x)) - z_g(L_i^{-1}(x))) \| \}\\
\le & s^*  \; \underset{x \in I_i} \sup \{ \| y_f(L_i^{-1}(x)) - y_g(L_i^{-1}(x)) +
 z_f(L_i^{-1}(x)) - z_g(L_i^{-1}(x)),\\
&  z_f(L_i^{-1}(x)) - z_g(L_i^{-1}(x)) \| \}\\
\le &  s^* \; \rho (f,g)
\end{split}
\end{equation*}
 where, in view of the conditions on $ \alpha_i, \beta_i, \gamma_i$  in Section 2.1,
$ s^* = \underset{1\le i \le N} \max
 \{ |\alpha_{i}|, |\beta_{i}|, |\gamma_{i}| \} < \; 1$. This shows that
 $T$ is a contraction mapping. By fixed point theorem, $T$ has a unique fixed point i.e. there exists
 a unique vector valued function $ f \in \mathcal{F}$ such that for all $ x \in I,$
$ (Tf)(x) = f(x).$ Now, for all $i = 1,2,\dots,N-1,$
\begin{equation*}
 f(x_{i}) = (Tf)(x_{i}) 
= F_{i+1} (L_{i+1}^{-1}(x_{i}), y(L_{i+1}^{-1}(x_{i}), z(L_{i+1}^{-1}(x_{i}))
= F_{i+1} (x_0,y_0,z_0) = (y_i,z_i),
\end{equation*}
which establishes that $f$ is the function interpolating the data
$\{(x_i,y_i,z_i) \; | \; i = 0,1,\dots,N \}$.
\par
It remains to show that the graph $\tilde{G}$ of the vector valued function $f$ is the
attractor of the IFS defined in (\ref{chf4}). To this end, observe that for all $x \in I,
i=1,2,\dots,N, $  and $f \in \mathcal{F},$
\begin{equation*}
(Tf)(L_i(x)) = F_{i} (x,y,z) =(\alpha_i y+ \beta_i z+c_i x+d_i,\; \gamma_i z+e_i x+f_i)
\end{equation*}
and 
\begin{equation*}
 \omega_{i} ( x,  y, z)  = (L_i(x) ,  F_{i}(x,y,z)) =  (L_i(x) , Tf(L_i(x)))
                    =  (L_i(x) , f(L_i(x)))
\end{equation*}
which implies that $\tilde{G}$ satisfies the invariance property, i.e.
$ \tilde{G} = \underset{i=1}{\overset{N}\bigcup} \omega_{i}(\tilde{G}).$
Since the  nonempty compact set that  satisfies the invariance property 
is unique, it follows that $ G = \tilde{G}.$ This proves 
$G$ is the graph of the vector valued function $f$ such that
$ G = \{ (x,y,z) | x \in I \}.$
\end{proof}
 Let the vector valued function $ f : I \rightarrow D $
 in Proposition \ref{prop2} be written as $ f(x)=(f_1(x), f_2(x))$.
 The required CHFIF is now  defined as follows:
\begin{Def}
Let $\{(x,f_1(x)) : x \in I \}$ be the projection of the attractor $G$ (c.f. (\ref{attractor})) 
on $\mathbb{R}^2$. Then, the function $f_1(x)$  is called
{\it { coalescence hidden variable FIF}} (CHFIF) for the given interpolation data
$\{ (x_i,y_i) \;| i=0,1,\dots,N \}.$ \end{Def}
\begin{Rem}
1. Although, the attractor $G$ (c.f. (\ref{attractor})) of the
IFS  defined in (\ref{chf4}) is a union of affine transformations of itself,
the projection of the attractor is not always union of affine
transformations of itself. Hence, CHFIFs are generally
non-self-affine in nature.\\
2. By choosing  $y_i = z_i$ and $\alpha_i+\beta_i =\gamma_i,$ 
CHFIF $f_1(x)$  obtained as the projection on
$\mathbb{R}^2$ of the attractor of the IFS (\ref{chf4})  
coincides with a self-affine FIF $f_2(x)$  for the same interpolation data.
Hence, the CHFIF is self-affine in this
case, in contrast to the observation of Barnsley\cite{B1} that the graph
of a HFIF is not self-similar or self-affine or self-anything.\\
3. For a given set of interpolation data, if an extra dimension is added to construct
the CHFIF, we have 1 free variable in the $3^{rd}$ co-ordinate
whereas, in the $2^{nd}$ co-ordinate, we have  1 free variable and 1 constrained variable.
In the case $y_i = z_i$, the resulting scaling factor of CHFIF is  $\alpha_i+\beta_i$.
As $|\alpha_i| < 1$ and $|\beta_i| + |\gamma_i| < 1$, taking
$|\gamma_i| < \epsilon$ for sufficiently small $\epsilon$,
the scaling factor of the CHFIF is found to lie
between $-2^+$ and $2^-$.\\
4. If extra $n$ dimensions  are added to interpolation data to get the CHFIF,
$(n+2)^{th}$ co-ordinate has 1 free variable,
$(n+1)^{th}$ co-ordinate has 1 free variable and may have at most 1 constrained
free variable, $n^{th}$ co-ordinate has 1 free variable and may have at most 2 constrained
free variables, $\dots.$ Continuing, the $2^{nd}$ co-ordinate has 1 free variable
and may have at most  $n$ constrained free variables.
So, in this extension, there are $n$ free variables and at most $(1+2+3+\dots+n)$ free
variables in the CHFIF. Due to the restrictions on free variables
and constrained free variables, the scaling factor of the CHFIF lies
between $-(n+1)^+$ and $(n+1)^-$. Thus, one can expect a wider range of
CHFIFs in higher dimension extensions.
\end{Rem}   
\section{ SMOOTHNESS ANALYSIS OF CHFIF}
\setcounter{equation}{0}
In this section, the smoothness of CHFIFs are studied by using their operator 
approximations. The H\"{o}lder exponent of  CHFIFs are calculated
in the proof of our main Theorems \ref{th1}-\ref{th3}.\\
We take the interpolation data on X-axis as $0 = x_0 < x_1 < \dots < x_N = 1.$
Let the function $F_i$ of the IFS (\ref{chf4}) be of the form
\begin{equation}\label{chf7}
F_i(x,y,z) = (\alpha_iy+\beta_iz+p_i(x),\; \gamma_{i}z+q_i(x))
\end{equation}
where $|\alpha_i| < 1, \; |\beta_i| + |\gamma_i| < 1, p_i \in Lip\lambda_i
( 0 < \lambda_i \le 1)$ and $ q_i \in Lip\mu_i (0 < \mu_i \le 1).$
From (\ref{chf6}) and (\ref{chf7}), for $ x \in I_i$, the fixed point f of T satisfies 
\begin{equation*}
\begin{split}
& Tf(x) = F_i(L_i^{-1}(x), f_1(L_i^{-1}(x)), f_2(L_i^{-1}(x)))\\
\Rightarrow & f(x) = F_i(L_i^{-1}(x), f_1(L_i^{-1}(x)), f_2(L_i^{-1}(x)))\\
\Rightarrow & (f_1(x),f_2(x)) = ( \alpha_i f_1 (L_i^{-1}(x)) +
\beta_i f_2 (L_i^{-1}(x)) + p_i
(L_i^{-1}(x)), \gamma_i f_2 (L_i^{-1}(x)) + q_i (L_i^{-1}(x))
\end{split}
\end{equation*}
Consequently, for all $ x \in I,$
\begin{equation*}
(f_1(L_i(x)), f_2(L_i(x))) = (\alpha_{i}f_1(x)+\beta_{i}f_2(x)
+p_i(x),\; \gamma_{i}f_2(x) + q_i(x))
\end{equation*}
Following Proposition \ref{prop2}, the CHFIF in this
case can be written as 
\begin{equation}\label{chf8}
f_1(L_i(x)) = \alpha_{i}f_1(x)+\beta_{i}f_2(x)+p_i(x)
\end{equation}
where, the self-affine fractal function $f_2(x)$  is given by 
\begin{equation}\label{chf9}
f_2(L_i(x)) = \gamma_{i}f_2(x) + q_i(x)
\end{equation}
Let $I_{r_1} = [x_{{r_1}-1}, x_{r_1}] = L_{r_1}(I).$ Then, 
$I_{r_1} = L_{r_1}(0) + |I_{r_1}| I,$ where $ |x_{r_1} - x_{{r_1}-1}|$ 
is the length of $ I_{r_1}, \; 1 \le r_1 \le N.$ Similarly,
$I_{{r_1}{r_2}} = L_{r_2}(0) + |I_{r_2}| L_{r_1}(I)= L_{r_2} \circ L_{r_1}(I)
 = L_{{r_1}{r_2}}(I),$
 where  $ |I_{{r_1}{r_2}}| = |I_{r_1}| . |I_{r_2}|, \; 1 \le r_1,r_2 \le N.$ 
In general,
\begin{equation}\label{chf10}
I_{{r_1}{r_2}\dots {r_m}} = L_{r_m}(0) + |I_{r_m}|
I_{{r_1}{r_2}\dots {r_{m-1}}} =\\ L_{r_m} \circ L_{r_{m-1}} \circ \dots
\circ L_{r_1}(I) = L_{{r_1}{r_2}\dots {r_m}}
\end{equation}
where, $ |I_{{r_1}{r_2}\dots {r_m}}| = |I_{r_1}| . |I_{r_2}| \dots |I_{r_m}| $ and
$ 1 \le r_1, r_2, \dots, r_m \le N.$\\
\par
 We need the following lemmas for our main results:
\begin{lm}\label{lm1}
Let  $f_1$ be defined as in (\ref{chf8}) and
$b_{{r_1}{r_2}\dots {r_m}} = \int_{I_{{r_1}{r_2}\dots {r_m}}} f_1(x)dx.$
Then,
\begin{equation}\label{chf11}
b_{{r_1}{r_2}\dots {r_m}}= \underset{k=1}{\overset{m}\sum} \underset{j=k+1}{\overset{m}\prod}
(|I_{r_j}|\alpha_{r_j})|I_{r_k}|(\int_{I_{{r_1}{r_2}\dots {r_{k-1}}}} p_{r_k}(\xi) d\xi +
 \beta_{r_k} a_{{r_1}{r_2}\dots {r_m}})+
\underset{j=1}{\overset{m}\prod} (|I_{r_j}|\alpha_{r_j})\int_0^{1} f_1(\xi) d\xi
\end{equation}
where, $ I_{r_0} = I \; \text{and} \;  a_{{r_1}{r_2}\dots {r_m}} = \int_{I_{{r_1}{r_2}\dots 
{r_m}}} f_2(x)dx.$
\end{lm}
\begin{proof} Since, $b_{{r_1}{r_2}\dots {r_m}} = \int_{L_{r_m}(0)+I_{{r_1}{r_2}\dots {r_{m-1}}}} 
f_1(x)dx$, a change of the variable $x$ by  $ x = L_{r_m}(0) + |I_{r_m}| \xi$ gives
\begin{equation*}
\begin{split}
b_{{r_1}{r_2}\dots {r_m}} &= \int_{I_{{r_1}{r_2}\dots {r_{m-1}}}}
f_1( L_{r_m}(0) + |I_{r_m}| \xi) |I_{r_m}| d\xi
=|I_{r_m}| \int_{I_{{r_1}{r_2}\dots {r_{m-1}}}} f_1( L_{r_m}(\xi))d\xi\\
&= |I_{r_m}| \int_{I_{{r_1}{r_2}\dots {r_{m-1}}}} (\alpha_{r_m}f_1(\xi) +
 \beta_{r_m}f_2(\xi) + p_{r_m}(\xi)) d\xi \\
&= |I_{r_m}| [ \int_{I_{{r_1}{r_2}\dots {r_{m-1}}}} p_{r_m}(\xi) d\xi + \beta_{r_m}
 \int_{I_{{r_1}{r_2}\dots {r_{m-1}}}} f_2(\xi) d\xi + \alpha_{r_m}
 \int_{I_{{r_1}{r_2}\dots {r_{m-1}}}} f_1(\xi) d\xi ]\\
&=|I_{r_m}| [\int_{I_{{r_1}{r_2}\dots {r_{m-1}}}} p_{r_m}(\xi) d\xi + \beta_{r_m}
a_{{r_1}{r_2}\dots {r_{m-1}}}] + |I_{r_m}| |I_{r_{m-1}}|\alpha_{r_m}
[\int_{I_{{r_1}{r_2}\dots {r_{m-2}}}} p_{r_{m-1}}(\xi) d\xi \\
&+ \beta_{r_{m-1}} a_{{r_1}{r_2}\dots {r_{m-2}}} + \alpha_{r_{m-1}}
\int_{I_{{r_1}{r_2}\dots {r_{m-2}}}} f_1(\xi) d\xi] = \dots \\
&=\underset{k=1}{\overset{m}\sum} \underset{j=k+1}{\overset{m}\prod}
(|I_{r_j}|\alpha_{r_j})|I_{r_k}|(\int_{I_{{r_1}{r_2}\dots {r_{k-1}}}} p_{r_k}(\xi) d\xi +
 \beta_{r_k} a_{{r_1}{r_2}\dots {r_m}})+
 \underset{j=1}{\overset{m}\prod} (|I_{r_j}|\alpha_{r_j})\int_0^{1} f_1(\xi) d\xi
\end{split}
\end{equation*}
\end{proof}
Since $f_1(x)$ is continuous, the integral average $b_{{r_1}{r_2}\dots {r_m}}/
|I_{{r_1}{r_2}\dots {r_m}}|$ can be taken as a good approximation of $f_1(x)$ in the
subinterval $I_{{r_1}{r_2}\dots {r_m}}$, when m is very large, leading to the
following definition of the approximating operator $\mathcal{Q}_m$ on the interval $I$:
\begin{Def}\label{def32} Let
\begin{equation}\label{chf12}
\mathcal{Q}_m(f_1, x) = \underset{{r_1}{r_2}\dots {r_m}=1}{\overset{N}\sum}
\chi_{I_{{r_1}{r_2}\dots {r_m}}}(x) \frac{b_{{r_1}{r_2}\dots {r_m}}}{|I_{{r_1}{r_2}\dots {r_m}}|}
\end{equation}
where, $I_{{r_1}{r_2}\dots {r_m}}$ is defined by (\ref{chf10}), $b_{{r_1}{r_2}\dots {r_m}}$
 is defined by (\ref{chf11}) and 
\begin{equation*}
\chi_{I_{{r_1}{r_2}\dots {r_m}}}(x)= \left\{\begin{array}{rlll}
1 &\quad x \in I_{{r_1}{r_2}\dots {r_m}}, \\
0 &\quad x \in I \setminus I_{{r_1}{r_2}\dots {r_m}}.
\end{array}\right.
\end{equation*}
\end{Def}
\begin{lm}\label{lm2}
The operator $\mathcal{Q}_m (f_1,x) $, given by (\ref{chf12}),  converges to $f_1(x)$ 
uniformly on I as \\ $m \rightarrow \infty.$\end{lm}
\begin{proof} The proof follows immediately by using Mean Value Theorem.\end{proof}
The following notations are needed throughout in the sequel:
$ \alpha = \max \{ |\alpha_i| : i = 1,2,\dots,N \},$
$\beta =\max \{ |\beta_i| : i = 1,2,\dots,N \}, \;
\gamma = \max \{ |\gamma_i| : i = 1,2,\dots,N \}, \;
\lambda = \min \{ \lambda_i : i = 1,2,\dots,N \},\;
\Omega_i =  \frac{|\alpha_i|}{|I_i|^{\lambda}},\;
\Omega = \max \{ \Omega_i : i = 1,2,\dots,N \}, \;
\mu = \min \{ \mu_i : i = 1,2,\dots,N \}, \;
\Gamma_i =  \frac{|\gamma_i|}{|I_i|^{\mu}}, \;
\Gamma = \max \{ \Gamma_i : i = 1,2,\dots,N \}, \; 
\Theta_i = \frac{|\alpha_i|}{|I_i|^{\mu}}, \;
\Theta = \max \{ \Theta_i : i = 1,2,\dots,N \},\; 
I_{\min} = \min \{ I_i : i = 1,2,\dots,N \}, \;   
I_{\max} = \max \{ I_i : i = 1,2,\dots,N \} $ and 
$\omega (f_1,t)$ = Modulus of continuity of $f_1(x).$
\par
Using  the above lemmas and notations, we now prove our  smoothness results 
according to the magnitude of $\Theta.$
\begin{thm}\label{th1} Let $f_1(x)$ be the CHFIF defined  by (\ref{chf8})
with $ \Theta < 1.$ Then, (a) for $\Omega \ne 1$ and $\Gamma \ne 1$, $f_1 \in Lip\delta$ 
(b) for $\Omega = 1$ or $\Gamma = 1,$  $\omega (f_1,t) = \bigcirc ( |t|^{\delta} \log|t|)$,
for suitable values of $\delta \in (0,1]$. 
\end{thm}
\begin{proof}
In order to calculate the H\"{o}lder exponent of  CHFIF 
$f_1$,  a suitable upper bound on the difference between $f_1(x)$ and
$f_1(\bar{x})$ for $ x, \bar{x} \in [0,1]$ is needed to be found.
 In view of Lemma \ref{lm2}, it is sufficient
to find an upper bound on  the difference between  functional values of  their operator
approximations $\mathcal{Q}_m (f_1, x) $   and $ \mathcal{Q}_m (f_1,\bar{x})$.
\par
 For $0\le x < \bar{x} \le 1,$ there exists a  least $m$ 
such that $ I_{{r_1}{r_2}\dots {r_m}}$ is the largest
interval contained in $[x, \bar{x} ].$
So, either $x$ or $\bar{x}  \in I_{{r_2}{r_3}\dots {r_m}}.$ Assume that,
$ x \in I_{s{r_2}\dots {r_m}},\; s \le r_1-1, \bar{x} \in I_{t{r_2}\dots {r_m}},
\; t \ge r_1+1, \; \text{or} \; \bar{x} \in I_{t^{\prime}{r_2+1}\dots {r_m}}, \;
1 \le t^{\prime} \le N.$ Let $n, m \in \N$ and $n > m$. Taking further
refinement of the above two intervals, we assume that
$ x \in I_{{u_1}\dots{u_{n-m}}s{r_2}\dots {r_m}}, \; \bar{x} \in
I_{{v_1}\dots{v_{n-m}}t{r_2^{\prime}}{r_3}\dots {r_m}}. $
It now follows that
\begin{equation*}
\begin{split}
\mathcal{Q}_n(f_1, x) &= \frac{1}{|I_{{u_1}\dots{u_{n-m}}s{r_2}\dots {r_m}}|}
 \int_{I_{{u_1}\dots{u_{n-m}}s{r_2}\dots {r_m}}} f_1(\xi)d\xi \\
  =& \frac{1}{{|I_{u_1}||I_{u_2}|\dots|I_{u_{n-m}}||I_s||I_{r_2}|\dots|I_{r_m}|}}
  [\underset{k=3}{\overset{m}\sum} \underset{j=k+1}{\overset{m}\prod}
  (|I_{r_j}|\alpha_{r_j})|I_{r_k}| \cdot \\
 & (\int_{I_{{u_1}\dots{u_{n-m}}s{r_2}\dots {r_{k-1}}}} p_{r_k}(\xi) d\xi +
   \beta_{r_k} a_{{u_1}\dots{u_{n-m}}s{r_2}\dots {r_{k-1}}})+
    \underset{j=3}{\overset{m}\prod} (|I_{r_j}|\alpha_{r_j})
    \int_{I_{{u_1}\dots{u_{n-m}}s{r_2}}} f_1(\xi) d\xi]
\end{split}
\end{equation*}
\begin{equation*}
\begin{split}
 =&\underset{k=3}{\overset{m}\sum} (\underset{j=k+1}{\overset{m}\prod}
 \alpha_{r_j}) \frac{1}{|I_{{u_1}\dots{u_{n-m}}s{r_2}\dots {r_{k-1}}}|}
  (\int_{I_{{u_1}\dots{u_{n-m}}s{r_2}\dots {r_{k-1}}}} p_{r_k}(\xi) d\xi +
  \beta_{r_k} a_{{u_1}\dots{u_{n-m}}s{r_2}\dots {r_{k-1}}})\\
 &+(\underset{j=3}{\overset{m}\prod} \alpha_{r_j}) \frac{1}
 {|I_{{u_1}\dots{u_{n-m}}s{r_2}}|}
 \int_{I_{{u_1}\dots{u_{n-m}}s{r_2}}} f_1(\xi) d\xi.
\end{split}
\end{equation*}
Similarly, the expression for $\mathcal {Q}_n (f_1,\bar{x})$ can be written as
\begin{equation*}
\begin{split}
\mathcal {Q}_n(f_1,\bar{x}) &= \underset{k=3}{\overset{m}\sum}
(\underset{j=k+1}{\overset{m}\prod} \alpha_{r_j})
\frac{1}{|I_{{v_1}\dots{v_{n-m}}t{r_2^{\prime}}\dots {r_{k-1}}}|}
(\int_{I_{{v_1}\dots{v_{n-m}}t{r_2^{\prime}}\dots {r_{k-1}}}} p_{r_k}(\xi) d\xi +
\beta_{r_k} a_{{v_1}\dots{v_{n-m}}t{r_2^{\prime}}\dots {r_{k-1}}})\\
&+(\underset{j=3}{\overset{m}\prod} \alpha_{r_j}) \frac{1}
{|I_{{v_1}\dots{v_{n-m}}t{r_2^{\prime}}}|}
\int_{I_{{v_1}\dots{v_{n-m}}t{r_2^{\prime}}}} f_1(\xi) d\xi
\end{split}
\end{equation*}
To estimate $| \mathcal {Q}_n(f_1,x) - \mathcal {Q}_n(f_1,\bar{x})|,$ observe that
\begin{equation}\label{chf13}
\begin{split}
 & \mathcal {Q}_n(f_1,x)-\mathcal {Q}_n(f_1,\bar{x}) \\ 
=& \underset{k=3}{\overset{m}\sum} (\underset{j=k+1}{\overset{m}\prod} \alpha_{r_j})
[\int_{I_{{u_1}\dots{u_{n-m}}s{r_2}\dots {r_{k-1}}}} \frac{p_{r_k}(\xi) d\xi}
{|I_{{u_1}\dots{u_{n-m}}s{r_2}\dots {r_{k-1}}}|} -
\int_{I_{{v_1}\dots{v_{n-m}}t{r_2^{\prime}}\dots {r_{k-1}}}}\frac{p_{r_k}(\xi) d\xi}
{|I_{{v_1}\dots{v_{n-m}}t{r_2^{\prime}} \dots {r_{k-1}}}|}] \\
&+(\underset{j=3}{\overset{m}\prod} \alpha_{r_j})
[\int_{I_{{u_1}\dots{u_{n-m}}s{r_2}}} \frac{f_1(\xi) d\xi}
{|I_{{u_1}\dots{u_{n-m}}s{r_2}}|} - \int_{I_{{v_1}\dots{v_{n-m}}t{r_2^{\prime}}}}
\frac{f_1(\xi) d\xi}{|I_{{v_1}\dots{v_{n-m}}t{r_2^{\prime}}}|}]\\
&+\underset{k=3}{\overset{m}\sum} (\underset{j=k+1}{\overset{m}\prod} \alpha_{r_j})
\beta_{r_k}[\frac {a_{{u_1}\dots{u_{n-m}}s{r_2}\dots {r_{k-1}}}}
{|I_{{u_1}\dots{u_{n-m}}s{r_2}\dots {r_{k-1}}}|} - \frac
{a_{{v_1}\dots{v_{n-m}}t{r_2^{\prime}}\dots {r_{k-1}}}}
{|I_{{v_1}\dots{v_{n-m}}t{r_2^{\prime}} \dots {r_{k-1}}}|}]
\end{split}
\end{equation}
Since  \cite{G},
\begin{equation*}
\begin{split}
a_{{r_1}{r_2}\dots {r_m}} & = \int_{I_{{r_1}{r_2}\dots {r_m}}} f_2(\xi)d\xi \\ 
 & = \underset{k=1}{\overset{m}\sum} \underset{j=k+1}{\overset{m}\prod}
(|I_{r_j}|\gamma_{r_j})|I_{r_k}|\int_{I_{{r_1}{r_2}\dots {r_{k-1}}}} q_{r_k}(\xi) d\xi +
\underset{j=1}{\overset{m}\prod} (|I_{r_j}|\gamma_{r_j})\int_0^{1} f_2(\xi) d\xi,
\end{split}
\end{equation*}
it follows that 
\begin{equation*}
\begin{split}
& \frac{a_{{u_1}\dots{u_{n-m}}s{r_2}\dots {r_{k-1}}}}
{|I_{{u_1}\dots{u_{n-m}}s{r_2}\dots {r_{k-1}}}|} =
\frac{1}{|I_{{u_1}\dots{u_{n-m}}s{r_2}\dots {r_{k-1}}}|}
\int_{I_{{u_1}\dots{u_{n-m}}{s}{r_2}\dots{r_{k-1}}}} f_2(\xi)d\xi\\
=&\frac{1}{|I_{{u_1}\dots{u_{n-m}}s{r_2}|}\dots|I_{r_3}|\dots|I_{r_{k-1}}|}
[\underset{l=3}{\overset{k-1}\sum} \underset{i=l+1}{\overset{k-1}\prod}
(|I_{r_i}|\gamma_{r_i})|I_{r_l}| \\
& \cdot \int_{I_{{u_1}\dots{u_{n-m}}{s}{r_2}\dots{r_{l-1}}}} q_{r_l}(\xi)d\xi +
\underset{i=3}{\overset{k-1}\prod}(|I_{r_i}|\gamma_{r_i})
\int_{I_{{u_1}\dots{u_{n-m}}{s}{r_2}}} f_2(\xi)d\xi] \\
=&\underset{l=3}{\overset{k-1}\sum} (\underset{i=l+1}{\overset{k-1}\prod}
\gamma_{r_i})\frac{1}{|I_{{u_1}\dots{u_{n-m}}{s}{r_2}\dots{r_{l-1}}}|}
\int_{I_{{u_1}\dots{u_{n-m}}{s}{r_2}\dots{r_{l-1}}}} q_{r_l}(\xi)d\xi \\
&+ (\underset{i=3}{\overset{k-1}\prod}\gamma_{r_i})
\frac{1}{|I_{{u_1}\dots{u_{n-m}}s{r_2}}|}
\int_{I_{{u_1}\dots{u_{n-m}}{s}{r_2}}} f_2(\xi)d\xi
\end{split}
\end{equation*}
Similarly,
\begin{equation*}
\begin{split}
\frac{a_{{v_1}\dots{v_{n-m}}t{r_2^{\prime}}\dots {r_{k-1}}}}
{|I_{{v_1}\dots{v_{n-m}}t{r_2^{\prime}}\dots {r_{k-1}}}|} =&
\underset{l=3}{\overset{k-1}\sum} (\underset{i=l+1}{\overset{k-1}\prod}
\gamma_{r_i})\frac{1}{|I_{{v_1}\dots{v_{n-m}}{t}{r_2^{\prime}}\dots{r_{l-1}}}|}
\int_{I_{{v_1}\dots{v_{n-m}}{t}{r_2^{\prime}}\dots{r_{l-1}}}} q_{r_l}(\xi)d\xi \\
&+ (\underset{i=3}{\overset{k-1}\prod}\gamma_{r_i})
\frac{1}{|I_{{v_1}\dots{v_{n-m}}t{r_2^{\prime}}}|}
\int_{I_{{v_1}\dots{v_{n-m}}{t}{r_2^{\prime}}}} f_2(\xi)d\xi
\end{split}
\end{equation*}
Consequently,
\begin{equation}\label{chf15}
\begin{split}
& |\frac{a_{{u_1}\dots{u_{n-m}}s{r_2}\dots {r_{k-1}}}}
{|I_{{u_1}\dots{u_{n-m}}s{r_2}\dots {r_{k-1}}}|} -
\frac{a_{{v_1}\dots{v_{n-m}}t{r_2^{\prime}}\dots {r_{k-1}}}}
{|I_{{v_1}\dots{v_{n-m}}t{r_2^{\prime}}\dots {r_{k-1}}}|} | \\
=&| \underset{l=3}{\overset{k-1}\sum} (\underset{i=l+1}{\overset{k-1}\prod} \gamma_{r_i})
[\int_{I_{{u_1}\dots{u_{n-m}}s{r_2}\dots {r_{l-1}}}} \frac{q_{r_l}(\xi) d\xi}
{|I_{{u_1}\dots{u_{n-m}}s{r_2}\dots {r_{l-1}}}|} -
\int_{I_{{v_1}\dots{v_{n-m}}t{r_2^{\prime}}\dots {r_{l-1}}}}\frac{q_{r_l}(\xi) d\xi}
{|I_{{v_1}\dots{v_{n-m}}t{r_2^{\prime}} \dots {r_{l-1}}}|}] \\
&+(\underset{i=3}{\overset{k-1}\prod} \gamma_{r_i})
[\int_{I_{{u_1}\dots{u_{n-m}}s{r_2}}} \frac{f_2(\xi) d\xi}
{|I_{{u_1}\dots{u_{n-m}}s{r_2}}|} - \int_{I_{{v_1}\dots{v_{n-m}}t{r_2^{\prime}}}}
\frac{f_2(\xi) d\xi}{|I_{{v_1}\dots{v_{n-m}}t{r_2^{\prime}}}|}]|\\
\le &  \underset{l=3}{\overset{k-1}\sum} \underset{i=l+1}{\overset{k-1}\prod} |\gamma_{r_i}|
[|\int_{I_{{u_1}\dots{u_{n-m}}s{r_2}\dots{r_{l-1}}}}
\frac{q_{r_l}(\xi)- q_{r_l}(x_{{r_3}\dots{r_{l-1}}})}
{|I_{{u_1}\dots{u_{n-m}}s{r_2}\dots {r_{l-1}}}|} d\xi|\\
& + |\int_{I_{{v_1}\dots{v_{n-m}}t{r_2^{\prime}}\dots{r_{l-1}}}}
\frac{q_{r_l}(\xi)- q_{r_l}(x_{{r_3}\dots{r_{l-1}}})}
{|I_{{v_1}\dots{v_{n-m}}t{r_2^{\prime}} \dots {r_{l-1}}}|} d\xi|] +
(\underset{i=3}{\overset{k-1}\prod} \gamma_{r_i})\cdot 2||f_2||_{\infty}\\
& \le  \underset{l=3}{\overset{k-1}\sum} (\underset{i=l+1}{\overset{k-1}\prod} |\gamma_{r_i}|)
 \; M_1 \; |I_{{r_3}\dots{r_{l-1}}}|^{\mu_{r_l}} + M_2 \;
(\underset{i=3}{\overset{k-1}\prod} \gamma_{r_i})
\end{split}
\end{equation}
where, $M_1$ is  Lipschitz bound and $M_2 = 2 \| f_2\|_{\infty}.$ 
Using (\ref{chf15}) in (\ref{chf13}),
\begin{equation*}
\begin{split}
&|\mathcal{Q}_n (f_1,x) - \mathcal{Q}_n (f_1,\bar{x})|
\le \underset{k=3}{\overset{m}\sum} (\underset{j=k+1}{\overset{m}\prod} |\alpha_{r_j}|)
[|\int_{I_{{u_1}\dots{u_{n-m}}s{r_2}\dots {r_{k-1}}}} \frac{p_{r_k}(\xi)
-p_{r_k} (x_{{r_3}\dots{r_{k-1}}})}
{|I_{{u_1}\dots{u_{n-m}}s{r_2}\dots {r_{k-1}}}|} d\xi|\\
&+|\int_{I_{{v_1}\dots{v_{n-m}}t{r_2^{\prime}}\dots {r_{k-1}}}}\frac{p_{r_k}(\xi)
- p_{r_k} (x_{{r_3}\dots{r_{k-1}}})}
{|I_{{v_1}\dots{v_{n-m}}t{r_2^{\prime}} \dots {r_{k-1}}}|} d\xi|]
+(\underset{j=3}{\overset{m}\prod} |\alpha_{r_j}|) \cdot 2||f_1||_{\infty}\\ &+
\underset{k=3}{\overset{m}\sum} (\underset{j=k+1}{\overset{m}\prod} |\alpha_{r_j}|)
|\beta_{r_k}|[\underset{l=3}{\overset{k-1}\sum} (\underset{i=l+1}{\overset{k-1}\prod}
|\gamma_{r_i}|)\; M_1\; |I_{{r_3}\dots{r_{l-1}}}|^{\mu_{r_l}} + M_2 \;
(\underset{i=3}{\overset{k-1}\prod} |\gamma_{r_i}|)]
\end{split}
\end{equation*}
\begin{equation*}
\begin{split}
\le & \underset{k=3}{\overset{m}\sum} (\underset{j=k+1}{\overset{m}\prod} |\alpha_{r_j}|)
\; M_3 \; |I_{{r_3}\dots{r_{k-1}}}|^{\lambda_{r_k}} + M_4 \;
(\underset{j=3}{\overset{m}\prod} |\alpha_{r_j}|) +
\underset{k=3}{\overset{m}\sum} (\underset{j=k+1}{\overset{m}\prod} |\alpha_{r_j}|)
|\beta_{r_k}| \cdot \\
 \quad & [\underset{l=3}{\overset{k-1}\sum} (\underset{i=l+1}{\overset{k-1}\prod}
|\gamma_{r_i}|)\; M_1\; |I_{{r_3}\dots{r_{l-1}}}|^{\mu_{r_l}} + M_2 \;
  (\underset{i=3}{\overset{k-1}\prod} |\gamma_{r_i}|)]
\end{split}
\end{equation*}
where,  $M_3$ is  Lipschitz bound and $ M_4 = 2 \| f_1 \|_{\infty}$.
From the above inequality it follows that
\begin{equation*}
\begin{split}
&|\mathcal{Q}_n (f_1,x) - \mathcal{Q}_n (f_1,\bar{x})|
\le   M_3 \;\underset{k=3}{\overset{m}\sum} (\underset{j=k+1}{\overset{m}\prod} |\alpha_{r_j}|)
 \frac {\underset{i^{\prime}=3}{\overset{m}\prod} |I_{r_{i^{\prime}}}|^{\lambda_{r_k}}}
 {\underset{j^{\prime}=k}{\overset{m}\prod} |I_{r_{j^{\prime}}}|^{\lambda_{r_k}}}
+M_4 \; (\underset{j=3}{\overset{m}\prod} |\alpha_{r_j}|) \\
& +\underset{k=3}{\overset{m}\sum} (\underset{j=k+1}{\overset{m}\prod} |\alpha_{r_j}|)
|\beta_{r_k}| \cdot 
[ M_1 \; \underset{l=3}{\overset{k-1}\sum} (\underset{i=l+1}{\overset{k-1}\prod}
|\gamma_{r_i}|) \quad \frac
{\underset{i^{\prime}=3}{\overset{m}\prod} |I_{r_{i^{\prime}}}|^{\mu_{r_l}}}
 {\underset{j^{\prime}=l}{\overset{m}\prod} |I_{r_{j^{\prime}}}|^{\mu_{r_l}}}
  + M_2 \; (\underset{i=3}{\overset{k-1}\prod} |\gamma_{r_i}|)]\\
\le & \frac{M_3}{({|I_{\min}|^{\lambda}})^3} \;(\underset{i^{\prime}=1}{\overset{m}\prod} 
 |I_{r_{i^{\prime}}}|^{\lambda})
\underset{k=3}{\overset{m}\sum} \underset{j=k+1}{\overset{m}\prod} \frac
{|\alpha_{r_j}|}{|I_{r_j}|^{\lambda}} + \frac{M_4}{({|I_{\min}|^{\lambda}})^2} 
\; (\underset{j=1}{\overset{m}\prod}|I_{r_j}|^{\lambda})
\underset{j=3}{\overset{m}\prod}\frac {|\alpha_{r_j}|}{|I_{r_j}|^{\lambda}} +
\underset{k=3}{\overset{m}\sum}|\beta_{r_k}|\underset{j=k+1}{\overset{m}\prod} |\alpha_{r_j}|
\cdot \\
&\left[\frac{M_1}{({|I_{\min}|^{\lambda}})^3} \; (\underset{i^{\prime}=1}{\overset{m}\prod} 
|I_{r_{i^{\prime}}}|^{\mu})
\underset{l=3}{\overset{k-1}\sum} \underset{i=l+1}{\overset{k-1}\prod}
\frac {|\gamma_{r_i}|}{|I_{r_i}|^{\mu}} + \frac{M_2}{({|I_{\min}|^{\lambda}})^2} \;
(\underset{i=1}{\overset{m}\prod}|I_{r_i}|^{\mu})
\underset{i=3}{\overset{k-1}\prod}\frac {|\gamma_{r_i}|}{|I_{r_i}|^{\mu}}
\right] \cdot  \underset{j=k}{\overset{m}\prod} \frac{1}{|I_{r_j}|^{\mu}}\\
\le & M_5 \; (\underset{i=1}{\overset{m}\prod} |I_{r_i}|^{\lambda})
[\underset{k=3}{\overset{m}\sum} \underset{j=k+1}{\overset{m}\prod} \Omega_{r_j} \; + \;
\underset{j=3}{\overset{m}\prod} \Omega_{r_j}] + M_6 \; \underset{k=3}
{\overset{m}\sum}|\beta_{r_k}|\underset{j=k+1}{\overset{m}\prod} 
\frac{|\alpha_{r_j}|}{|I_{r_j}|^{\mu}} \cdot \\
 \quad  & (\underset{i=1}{\overset{m}\prod}|I_{r_i}|^{\mu})
 \left[\underset{l=3}{\overset{k-1}\sum} \underset{i=l+1}{\overset{k-1}\prod} \Gamma_{r_i} \; + \;
 \underset{i=3}{\overset{k-1}\prod} \Gamma_{r_i}\right]
\end{split}
\end{equation*}
where, $M_5 = Max \{\frac{M_3}{({|I_{\min}|^{\lambda}})^3} \frac{M_4}{({|I_{\min}|^{\lambda}})^2}\}$
and $M_6 = Max \{\frac{M_1} {({|I_{\min}|^{\lambda}})^4},
\frac{M_2}{({|I_{\min}|^{\lambda}})^3}\}.$ The above inequality gives 
\begin{equation}\label{chf17}
|\mathcal{Q}_n (f_1,x) - \mathcal{Q}_n (f_1,\bar{x})|
\le  M_5 \; |x-\bar{x}|^{\lambda}(\underset{k=2}{\overset{m}\sum} \Omega^{m-k}) +
M_6 \; \beta |x-\bar{x}|^{\mu} \underset{k=3}{\overset{m}\sum} \Theta^{m-k}
(\underset{l=3}{\overset{k-1}\sum} \Gamma^{k-l}) 
\end{equation}
 Since $\Theta < 1,$ (\ref{chf17}) further reduces to
\begin{equation}\label{chf18}
|\mathcal{Q}_n (f_1,x) - \mathcal{Q}_n (f_1,\bar{x})|
\le  M_5\;  |x-\bar{x}|^{\lambda}(\underset{k=2}{\overset{m}\sum} \Omega^{m-k}) +
M_6 \; \frac{\beta}{1-\Theta}|x-\bar{x}|^{\mu}(\underset{k=1}{\overset{m-3}\sum} \Gamma^k)
\end{equation}
{\bf Case (a).} $\Omega \ne 1$ and $\Gamma \ne 1$: The desired H\"{o}lder exponents
are found individually for each of the following subcases\\
{\bf I.} $\Omega < 1$ and $\Gamma < 1 :$
$|\mathcal{Q}_n (f_1,x) - \mathcal{Q}_n (f_1,\bar{x})| \le
 \frac{M_5}{1-\Omega}|x-\bar{x}|^{\lambda} + \frac{M_6\;\beta}{(1-\Theta)(1-\Gamma)}
 |x-\bar{x}|^{\mu}$ $\le  M_7 |x-\bar{x}|^{\delta_1}$
where, $M_7 = \max\{ \frac{M_5}{1-\Omega}, \frac{M_6\;\beta}{(1-\Theta)(1-\Gamma)}\}$ and 
$\delta_1 = \min (\lambda, \mu) $. Thus, as  $n\rightarrow \infty$,
the above inequality together with Lemma \ref{lm2} gives  $f_1 \in Lip \delta$
with $\delta= \delta_1.$\\
{\bf II.} $\Omega > 1$ and $\Gamma > 1:$
\begin{equation}\label{chf19}
|\mathcal{Q}_n (f_1,x) - \mathcal{Q}_n (f_1,\bar{x})| \le
M_5 \;|x-\bar{x}|^{\lambda}\; m\; \Omega^m +\frac{M_6\; \beta}{1-\Theta}|x-\bar{x}|^{\mu}
 \;m \;\Gamma^m
\end{equation}
Suppose $\tau_1 > 0$ such that $|x-\bar{x}|^{\lambda}\; m\; \Omega^m \; \le
\; |x-\bar{x}|^{\tau_1}$. Then,\\
\begin{equation}\label{chf20}
\tau_1 \le \lambda + \frac{m \; \log{\Omega}}{\log{|x-\bar{x}|}}
\end{equation}
Further, $|I_{{r_1}\dots{r_m}}| \le |x-\bar{x}| < 1 \; \Rightarrow
|I_{\min}|^m \le |x-\bar{x}| \; \Rightarrow \frac{1}{m\log|I_{\min}|}
 \ge \frac{1}{\log|x-\bar{x}|}.$\\
Also, $ \Omega \le \frac{\alpha}{|I_{\min}|^{\lambda}} \;\Rightarrow
\log{\Omega} \le \log{\alpha} - \lambda \log|I_{\min}|.$  Therefore, by (\ref{chf20}),
$\tau_1 \le \frac{\log{\alpha}}{\log|I_{\min}|} $.\\
Similarly, if $\tau_2 > 0$ is such that $|x-\bar{x}|^{\mu}\; m\; \Gamma^m \; \le
\; |x-\bar{x}|^{\tau_2}$, then
$\tau_2 \le \frac{\log{\gamma}}{\log|I_{\min}|} $.\\
Let $ \tau_3 = \min \{ \frac{\log{\alpha}}{\log|I_{\min}|}, \; \;
\frac{\log{\gamma}}{\log|I_{\min}|} \}. $
From (\ref{chf19}), for any $\delta_2 \le \tau_3,$  
$|\mathcal{Q}_n (f_1,x) - \mathcal{Q}_n (f_1,\bar{x})| \le
M_8 \; |x-\bar{x}|^{\delta_2},$ where $M_8 = \max \{ M_5, \frac{M_6\; \beta}{1-\Theta}\}.$
 Now, the last inequality  together with 
Lemma \ref{lm2} gives  $ f_1 \in Lip \delta$ with $\delta= \delta_2.$\\
{\bf III.} $\Omega > 1$ and $\Gamma < 1:$
$|\mathcal{Q}_n (f_1,x) - \mathcal{Q}_n (f_1,\bar{x})| \le
M_5 |x-\bar{x}|^{\tau_1} + \frac{M_6\; \beta}
{(1-\Theta)(1-\Gamma)}|x-\bar{x}|^{\mu} \le M_{9} \; |x-\bar{x}|^{\delta_3}$
where, $M_9 = \max \{ M_5, \frac{M_6\; \beta}{(1-\Theta)(1-\Gamma)}\}$ 
and  $ \delta_3 = \min (\tau_1, \mu).$ Thus,
as  $n\rightarrow \infty$, the last inequality together with Lemma \ref{lm2} gives
 $ f_1 \in Lip \delta$ with $\delta= \delta_3.$\\
{\bf IV.} $\Omega < 1$ and $\Gamma > 1:$
$|\mathcal{Q}_n (f_1,x) - \mathcal{Q}_n (f_1,\bar{x})| \le
\frac{M_5}{1-\Omega} |x-\bar{x}|^{\lambda} + \frac{M_6\; \beta}
{1-\Theta}|x-\bar{x}|^{\tau_2}
\le M_{10} \; |x-\bar{x}|^{\delta_4} $
where, $ M_{10} = \max\{ \frac{M_5}{1-\Omega}, \frac{M_6\; \beta}
{1-\Theta} \}$ and  $ \delta_4 = \min ( \lambda , \tau_2).$ So, as  
$n\rightarrow \infty$, the above inequality together with Lemma \ref{lm2} gives
 $ f_1 \in Lip\delta$ with  $\delta= \delta_4$.\\
{\bf Case (b).} $\Omega = 1$ or $\Gamma = 1$: The desired H\"{o}lder exponents
are found individually for each of the following subcases\\
{\bf I.} $\Omega = 1$ and $\Gamma \le 1$ or $\Omega < 1$ and $\Gamma = 1:$
For $\Omega = 1$ and $\Gamma =1$, 
\begin{equation*}
\begin{split}
|\mathcal{Q}_n (f_1,x) - \mathcal{Q}_n (f_1,\bar{x})| &\le
(M_5|x-\bar{x}|^{\lambda} + \frac{M_6\;\beta}{1-\Theta}|x-\bar{x}|^{\mu}) \cdot (m-1)\\
 &\le (M_5|x-\bar{x}|^{\lambda} + \frac{M_6\;\beta}{1-\Theta}|x-\bar{x}|^{\mu})
\frac {\log|x-\bar{x}|}{\log|I_{\max}|}\\
&\le M_{11} \; (\log|x-\bar{x}|) |x-\bar{x}|^{\delta_1}
\end{split}
\end{equation*}
where, $M_{11} = \frac{M_8}{\log|I_{\max}|}$. As  $n\rightarrow \infty$, the 
last inequality together with Lemma \ref{lm2} gives
$\omega (f_1,t) =\bigcirc ( |t|^{\delta}\log|t|)$ with $\delta= \delta_1$.
For $\Omega = 1$ and $\Gamma < 1$, $|\mathcal{Q}_n (f_1,x) - \mathcal{Q}_n (f_1,\bar{x})| \le
\frac{M_5}{\log|I_{\max}|} |x-\bar{x}|^{\lambda} \log{|x-\bar{x}|} +
\frac{M_6\; \beta}{(1-\Theta)(1-\Gamma)}|x-\bar{x}|^{\mu}
\le M_{12} \; |x-\bar{x}|^{\delta_1} (1+ \log|x-\bar{x}|),$ 
where $M_{12} = \max \{ \frac{M_5}{\log|I_{\max}|}, \frac{M_6\; \beta}{(1-\Theta)(1-\Gamma)} \}.$
Hence, as  $n\rightarrow \infty$, the above inequality together with Lemma \ref{lm2} gives
 $\omega (f_1,t) = \bigcirc ( |t|^{\delta}(1+\log|t|))
\equiv \bigcirc ( |t|^{\delta}\log|t|)$ with $\delta= \delta_1$.
The estimate for $\Omega < 1$ and $\Gamma = 1$ follows using analogous arguments.\\
{\bf II.} $\Omega > 1$ and  $\Gamma = 1:$
$|\mathcal{Q}_n (f_1,x) - \mathcal{Q}_n (f_1,\bar{x})| \le
M_5 \; |x-\bar{x}|^{\tau_1} +
\frac{M_6 \; \beta}{(1-\Theta)\log|I_{\max}|} |x-\bar{x}|^{\mu} \log|x-\bar{x}| 
\le M_{13} |x-\bar{x}|^{\delta_3} (1+ \log|x-\bar{x}|),$ where $M_{13} = \max\{ M_5,
\frac{M_6 \; \beta}{(1-\Theta)\log|I_{\max}|} \}$.
Making  $n\rightarrow \infty$, the above inequality together with Lemma \ref{lm2} gives 
 $\omega (f_1,t) = \bigcirc ( |t|^{\delta}(1+\log|t|))
\equiv \bigcirc ( |t|^{\delta}\log|t|)$  with $\delta= \delta_3$.\\
{\bf III.} $\Omega = 1$ and $\Gamma > 1:$ 
$|\mathcal{Q}_n (f_1,x) - \mathcal{Q}_n (f_1,\bar{x})| \le
\frac{M_5}{\log|I_{\max}|} |x-\bar{x}|^{\lambda} \log{|x-\bar{x}|} +
\frac{M_6 \; \beta}{1-\Theta} |x-\bar{x}|^{\tau_2} 
\le M_{14} |x-\bar{x}|^{\delta_4} (1+ \log|x-\bar{x}|),$
where $M_{14} = \max \{ \frac{M_5}{\log|I_{\max}|}, \frac{M_6 \; \beta}{1-\Theta} \}$.
So, as  $n\rightarrow \infty$, the above inequality together with Lemma \ref{lm2} gives 
 $\omega (f_1,t) = \bigcirc ( |t|^{\delta}(1+\log|t|))
\equiv \bigcirc ( |t|^{\delta}\log|t|)$ with $\delta= \delta_4$.\\
Theorem \ref{th1} now follows from the above cases with suitable $\delta$
as found in various subcases.\end{proof}
 The smoothness results for the class of  CHFIFs when 
$\Theta =1$ are given by the  following:
\begin{thm}\label{th2} Let $f_1(x)$ be the  CHFIF defined  by (3.2)
with $ \Theta = 1.$  Then, (a) for $\Omega \ne 1$ and $\Gamma \ne 1$, 
$\omega (f_1,t) = \bigcirc ( |t|^{\delta}\log|t|)$ (b) for 
$\Omega = 1$ or  $\Gamma = 1$,
$\omega (f_1,t) = \bigcirc ( |t|^{\delta} (\log|t|)^2),$ for suitable values of $\delta \in (0,1]$.
\end{thm}
\begin{proof} Since  $\Theta =1 $,  (\ref{chf17}) gives,
\begin{equation}\label{chf21}
|\mathcal{Q}_n (f_1, x) - \mathcal{Q}_n (f_1,\bar{x})| 
\le M_5\;  |x-\bar{x}|^{\lambda}(\underset{k=2}{\overset{m}\sum} \Omega^{m-k}) +
\frac{M_6 \; \beta}{\log|I_{\max}|}|x-\bar{x}|^{\mu} \; \log|x-\bar{x}|
(\underset{k=1}{\overset{m-3}\sum} \Gamma^k)
\end{equation}
The rest of proof is similar to that of Theorem \ref{th1} with the respective values of
$\delta$ as in different cases of Theorem \ref{th1}.\end{proof}
Finally, the  smoothness results for the class of  CHFIFs for  $\Theta > 1 $ 
are given by the  following:
\begin{thm}\label{th3}
 Let $f_1(x)$ be the  CHFIF defined  by (\ref{chf8}) with $ \Theta > 1.$ 
Then, (a) for $\Omega \ne 1$ and $\Gamma \ne 1$, $f_1 \in Lip\delta$  
(b) for $\Omega = 1$ or $\Gamma = 1,$ 
$\omega (f_1,t) = \bigcirc ( |t|^{\delta} \log|t|),$ for suitable values of $\delta \in (0,1]$.
\end{thm}
\begin{proof} Inequality (\ref{chf17}) for  $\Theta > 1 $ gives
\begin{equation}\label{chf22}
|\mathcal{Q}_n (f_1, x) - \mathcal{Q}_n (f_1,\bar{x})|  \le
M_5\;  |x-\bar{x}|^{\lambda}(\underset{k=2}{\overset{m}\sum} \Omega^{m-k}) +
M_6 \; \beta |x-\bar{x}|^{\mu} \cdot m \Theta^m \cdot
(\underset{k=1}{\overset{m-3}\sum} \Gamma^k) 
\end{equation}
Let  $\tau_4 > 0$ be  such that $|x-\bar{x}|^{\mu}\; m\; \Theta^m \; \le
\; |x-\bar{x}|^{\tau_4}$. Then,
\begin{equation*}
\tau_4 \le \mu + \frac{m \; \log{\Theta}}{\log{|x-\bar{x}|}}
\le \frac{\log{\alpha}}{\log|I_{\min}|}.
\end{equation*}
Since  $\tau_1$ in Theorem \ref{th1} satisfies $\tau_1 \le \frac{\log{\alpha}}{\log|I_{\min}|}$,
we can choose  $ \tau_4 = \tau_1$ so that (\ref{chf22}) reduces to
\begin{equation}\label{chf23}
|\mathcal{Q}_n (f_1, x) - \mathcal{Q}_n (f_1,\bar{x})|  \le
M_5\;  |x-\bar{x}|^{\lambda}(\underset{k=2}{\overset{m}\sum} \Omega^{m-k}) +
M_6 \; \beta |x-\bar{x}|^{\tau_1} \; (\underset{k=1}{\overset{m-3}\sum} \Gamma^k)
\end{equation} 
The rest of the  proof is similar to that of Theorem \ref{th1}, by considering
(\ref{chf23}) in place of (\ref{chf18}). As in Theorem \ref{th1}, the value of $\delta$ in different cases
are given by Case (a):  I. $\delta = \delta_5 = \min (\lambda, \tau_1)$,
II. $\delta = \delta_6$ where,  $\delta_6 \le \frac{\log{\alpha}{\gamma}}
{\log|I_{\min}|} - \mu,$ III. $\delta = \delta_7$ where, $\delta_7 \le 
\frac{\log{\alpha}}{\log|I_{\min}|}$, IV. 
$\delta = \delta_8 = \min (\lambda, \delta_6),$ and  Case (b):
I. $\delta = \delta_5$, II. $\delta = \delta_7$, 
III. $\delta = \delta_8$.
\end{proof}
\begin{Rem}\label{re2}
1. It follows from Theorems \ref{th1}-\ref{th3}, that the smoothness of the CHFIF depends 
on the free variables $ \alpha_i, \gamma_i$ and  the Lipschitz exponents
$\lambda_i$ and $\mu_i$.\\
2. If $p_i(x)$ and $q_i(x)$ belong to the same function space, then
$\lambda_i = \mu_i \Rightarrow \lambda = \mu \Rightarrow \Omega_i = \Theta_i
\Rightarrow \Omega = \Theta.$
Thus, in this case, there are only
three subcases in each of Theorems \ref{th1}-\ref{th3}, depending on magnitude of $\Gamma.$
The  CHFIF $f_1(x)$ is not self-affine
if either $y_i \neq z_i$ for $i=0,1,2,\dots,N$ or 
$p_i(x) \neq q_i(x)$ or $\alpha_i + \beta_i \neq \gamma_i$ for
$ i = 1,2,\dots,N.$ Thus, we need to choose $\lambda \neq \mu,$
to obtain all the nine subcases of Theorems \ref{th1}-\ref{th3}.\\
3. Let $\lambda = \mu$ and $\Theta < 1$. Then, $ \Omega < 1.$
Theorem \ref{th1} now gives the following  smoothness results
depending on the magnitude of $\Gamma$ for  CHFIF $f_1(x).$
(A) For $\Gamma < 1, $  $f_1(x) \in \text{Lip}\mu$, since  $\delta_1 = \mu$ in this case.
(B) For $\Gamma = 1, $  $\omega(f_1, t) = \bigcirc (|t|^{\mu} \log|t|),$ since 
$\delta_1 = \mu$ in this case.
(C) For $\Gamma > 1, $  $f_1 \in \text{Lip}\tau_2 \; \text{where}
\; \tau_2 \le  \frac{\log \gamma}{log|I_{\min}|},$ since $ \delta_4 = \min (\lambda, \tau_2)
\le \tau_2$ in this case.\\
4. Suppose $\lambda = \mu$  and $\Theta = 1$. Then, $\Omega=1.$   Theorem \ref{th2} in this case
gives the smoothness result as follows:
(A) For $\Gamma < 1, $  $\omega(f_1, t) = \bigcirc (|t|^{\mu} (\log|t|)^2),$
since $\delta_1 = \mu$ in this case.
(B) For $\Gamma = 1, $  $\omega(f_1, t) = \bigcirc (|t|^{\mu} (\log|t|)^2),$ 
since $\delta_1 = \mu$ in this case.
(C) For $\Gamma > 1, $  $\omega(f_1, t) = \bigcirc (|t|^{\delta_4} (\log|t|)^2),$
where $ \delta_4 = \min (\lambda, \tau_2) \le \tau_2.$\\
5. Let $\lambda = \mu$  and $\Theta > 1$. Then, $ \Theta = 
\frac{\max\{|\alpha_i| : i=1,2,\dots,N\}}{|I_{\min}|^{\mu}} > 1$ which in turn implies 
$\frac{\log \alpha}{\log |I_{\min}|} < \mu.$ 
Since  $\delta_6 \le \frac{\log \gamma}{\log |I_{\min}|} +
(\frac{\log \alpha}{\log |I_{\min}|} - \mu) < \frac{\log \gamma}{\log |I_{\min}|}$
and $\tau_2 \le \frac{\log \gamma}{\log |I_{\min}|},$ we may choose $\delta_6 \le \tau_2.$
Further, $\tau_1 \le \frac{\log \alpha}{\log |I_{\min}|}$ implies  $\tau_1 < \mu.$
With  these inequalities,  the smoothness results as derived from Theorem \ref{th3}
in the case $\lambda = \mu$  and $\Theta= \Omega > 1$ are as follows:
(A) For $\Gamma < 1, $ $f_1 \in \text{Lip}\tau_1 \supseteq \text{Lip}\mu.$
(B) For $\Gamma = 1, $ $\omega(f_1, t) = \bigcirc (|t|^{\tau_1} \log|t|)$ 
which gives $\omega(f_1, t)= \bigcirc (|t|^{\mu} \log|t|).$
(C) For $\Gamma > 1, $ $f_1 \in \text{Lip}\delta_6 \supseteq \text{Lip}\tau_2.$\\
6. If $f_1(x)= f_2(x),$ then $f_1(x)$  is also self-affine and in such case, $y_i = z_i, \;
\alpha_i + \beta_i = \gamma_i$ and $p_i(x) = q_i(x).$ Hence, $ \lambda_i = \mu_i \Rightarrow
\lambda = \mu \Rightarrow \delta_1 = \mu$ and $ \Omega = \Theta.$
For self-affine function $f_1(x) = f_2(x),$  $f_1$ belongs to the  
intersection of the function spaces occurring for the same case((A), (B), or (C))
of Remarks 3-5 as above. We note that the intersection of these function spaces is independent
of $\Theta.$ Since, the class of  CHFIFs for  $\Theta < 1$ is
contained in the class of  CHFIFs for $\Theta = 1$ and $\Theta > 1$,
the smoothness results in \cite{G} for self-affine function $f_2(x)$ follows as
special case of our smoothness results derived in the above Remarks 3-5.
\end{Rem}
\section{FRACTAL DIMENSION AND  CHFIF }
\setcounter{equation}{0}
The following definitions are needed in the sequel:
 The conditions  $\Omega = 1$, $\Gamma = 1$ or $\Theta =1$  are called 
{\it {critical conditions}}.  The CHFIF $f_1(x)$  with any one of these condition
is called {\it {critical  CHFIF.}} Let $\mathcal{N}(A,\epsilon)$ be the smallest
number of closed balls of radius $\epsilon > 0,$  needed to
cover A. Then, the Fractal dimension of A is defined by
$ D_B(A) = \underset{\epsilon \rightarrow 0}{\lim}
 \frac{\log \mathcal{N}(A,\epsilon)}{- \log \epsilon}$,whenever the limit exists.
Our following theorems give bounds of the fractal dimension for the critical  CHFIFs.
\begin{thm}\label{th4} Let  CHFIF  $f_1(x)$ be defined by
(\ref{chf8}). Then, for the critical condition $\Omega = 1 $, 
\begin{equation}\label{fd}
1- \frac{\log \underset{k=1}{\overset{N}\sum} |\alpha_k|}{\log |I_{\max}|}
 \le D_B(graph(f_1)) \le 1 - \delta - \frac{\log N}{\log |I_{\max}|}
\end{equation}
and for the  critical condition $\Gamma = 1,$
\begin{equation}\label{fd1}
1- \frac{\log \underset{k=1}{\overset{N}\sum} |\gamma_k|}{\log |I_{\max}|}
 \le D_B(graph(f_1)) \le 1 - \delta - \frac{\log N}{\log |I_{\max}|}
\end{equation}
where, $\delta$ takes suitable values as in the subcases in Theorems \ref{th1}-\ref{th3}.  
\end{thm}
\begin{proof} Let $\Theta < 1$ and $\Omega = 1.$ 
Since $\omega (f_1,t) = \bigcirc ( |t|^{\delta_1}\log|t|),$ (c.f. Theorem \ref{th1}),
for all$ x \ne x^*$,  $ x, x^* \in I,$ there exist constants
 $C_1, C_2 $ such that
\begin{equation}\label{fd2}
C_1 |x-x^*|^{\delta_1} \le |f_1(x) - f_1(x^*)| \le C_2 |x-x^*|^{\delta_1}\log|x-x^*|
\end{equation}
Suppose, $ G_{r_1,r_2,\dots,r_m} = \{(x,f_1(x),f_2(x))\; |\; x \in I_{r_1,r_2,\dots,r_m} \}.$
Define, $|A|_X = \sup \{ |x-\bar{x}| \; | \; (x,y,z),$ $ (\bar{x},\bar{y},\bar{z}) \in A \},$
 $|A|_Y = \sup \{ |y-\bar{y}| \; | \; (x,y,z), (\bar{x},\bar{y},\bar{z}) \in A \},$ for
 any $ A \subset \mathbb{R}^3.$ Since, $|G_{r_1,r_2,\dots,r_m}|_X = |I_{r_1,r_2,\dots,r_m}|,$
 (\ref{fd2}) reduces to
\begin{equation}\label{fd3}
C_1 |I_{r_1,r_2,\dots,r_m}|^{\delta_1} \le |G_{r_1,r_2,\dots,r_m}|_Y \le
C_2 |I_{r_1,r_2,\dots,r_m}|^{\delta_1} \log |I_{r_1,r_2,\dots,r_m}|
\end{equation}
Choose $m$ large  such that $ |I_{\max}|^m < \frac{1}{2} \epsilon, \; \epsilon > 0.$
Since, $ \Omega_{r_j} \le \Omega = 1$ implies $ |\alpha_{r_j}| \le
|I_{r_j}|^{\lambda} \le |I_{r_j}|^{\delta_1}$ and $|I_{r_1,r_2,\dots,r_m}|
= |I_{r_1}| \cdot |I_{r_2}| \dots |I_{r_m}|,$ it follows by (\ref{fd3}) that
\begin{equation}\label{fd4}
 C_1 |\alpha_{r_1}|\cdot|\alpha_{r_2}|\dots |\alpha_{r_m}| \le |G_{r_1,r_2,\dots,r_m}|_Y
\le C_2 |I_{\max}|^{m\delta_1} \cdot m\log|I_{\max}|
\end{equation}
Taking summation over $r_1,r_2,\dots,r_m$ from $1$ to $N$ in 
(\ref{fd4}), 
\begin{equation*}
\begin{split}
\underset{r_1,r_2,\dots,r_m}{\sum} C_1
|\alpha_{r_1}|\cdot|\alpha_{r_2}|\dots |\alpha_{r_m}| |I_{\max}|^{-m} & \le
\underset{r_1,r_2,\dots,r_m}{\sum}|G_{r_1,r_2,\dots,r_m}|_Y |I_{\max}|^{-m} \\
 \le & \underset{r_1,r_2,\dots,r_m}{\sum} C_2 |I_{\max}|^{m(\delta_1-1)}
\cdot m\log|I_{\max}|
\end{split}
\end{equation*}
The above inequalities can be rewritten as
\begin{equation*}
\bar{C_1} |I_{\max}|^{-m} (|\alpha_1| + \dots |\alpha_N|)^m
\le \mathcal{N}(\text{graph}(f_1), \epsilon) \le
\bar{C_2} |I_{\max}|^{m(\delta_1-1)} \cdot m\log|I_{\max}| \cdot N^m
\end{equation*}
The inequalities (\ref{fd}) follow from the
last inequalities with $\delta = \delta_1.$ The proof of (\ref{fd1}) for
$\Theta < 1$, $\Gamma =1$ is  similar to the above case.
\par
Let $\Theta = 1$ and $\Omega= 1.$ Since $\omega (f_1,t) = 
\bigcirc ( |t|^{\delta_1}(\log|t|)^2)$ (c.f. Theorem \ref{th2}), 
for all $x \ne x^*$, $ x, x^* \in I,$ there exist constants
$C_3, C_4 $ such that
\begin{equation}\label{e1}
C_3 |x-x^*|^{\delta_1} \le |f_1(x) - f_1(x^*)| \le C_4 |x-x^*|^{\delta_1}(\log|x-x^*|)^2
\end{equation}
Now, using (\ref{e1}) in place of (\ref{fd2}) 
the above arguments give that there are constants $\bar{C_3}$ and $\bar{C_4}$ such that
\begin{equation*}
\bar{C_3} |I_{\max}|^{-m} (|\alpha_1| + \dots |\alpha_N|)^m
\le \mathcal{N}(\text{graph}(f_1), \epsilon) \le
\bar{C_4} |I_{\max}|^{m(\delta_1-1)} \cdot (m\log|I_{\max}|)^2 \cdot N^m
\end{equation*}
The proof of (\ref{fd})  for $\Theta = 1$ and $\Omega= 1$ follows
from the above inequalities  with $\delta= \delta_1$.
\par
The proof of  (\ref{fd})-(\ref{fd1}) is analogous in other cases.\end{proof}
\begin{thm}\label{th5}
Let  CHFIF $f_1(x)$ be defined in
(\ref{chf8}) with $\Theta = 1$. Then, for $\Omega \ne 1 $ or $\Gamma \ne 1$,
\begin{equation*}
1- \frac{\log \underset{k=1}{\overset{N}\sum} |\alpha_k|}{\log |I_{\max}|}
 \le D_B(graph(f_1)) \le 1 - \delta - \frac{\log N}{\log |I_{\max}|}
\end{equation*}
where, $\delta$ takes suitable values as in Theorem \ref{th2}.
\end{thm}
\begin{proof} The proof is similar to the case $\Theta < 1 $ of Theorem \ref{th4}.\end{proof}
Theorems \ref{th4}-\ref{th5} lead to the following bounds on
fractal dimension of equally spaced critical  CHFIFs.
\begin{corl}\label{cor1}
Let  CHFIF $f_1(x)$ be defined by (\ref{chf8}). 
Then, for  $\Theta =1$ or $\Omega =1$,
\begin{equation}\label{e2}
1+ \frac{\log \underset{k=1}{\overset{N}\sum} |\alpha_k|}{\log N}
 \le D_B(graph(f_1)) \le 2 - \delta.
\end{equation}
Further, for   $\Gamma = 1,$
\begin{equation}\label{e3}
1+ \frac{\log \underset{k=1}{\overset{N}\sum} |\gamma_k|}{\log N}
 \le D_B(graph(f_1)) \le 2 - \delta
\end{equation}
where, $\delta$ takes suitable values  as in Theorems \ref{th1}-\ref{th3}.
\end{corl}
\begin{corl}\label{cor2}
Let  the equidistant CHFIF  $f_1(x)$ be defined by (\ref{chf8}). Then,
$ D_B(graph(f_1)) = 1$ in the following cases:\\
1. $\Theta \le  1$, either $\delta = \delta_1 = 1$ or $\delta = \delta_3 = 1$ and either 
$\underset{k=1}{\overset{N}\sum} \; |\alpha_k|
 \le \; 1 $   or  $\underset{k=1}{\overset{N}\sum} \;|\gamma_k| \le 1$.\\ 
2. $\Theta \le  1$,  $\delta = \delta_4 = 1,$ and $\underset{k=1}{\overset{N}\sum}
|\alpha_k| \le 1.$\\
3. $\Theta = 1$,  $\delta = \delta_2 = \tau_3 = 1,$ 
and  $\underset{k=1}{\overset{N}\sum} \; |\alpha_k| \le 1 $.\\
4. $\Theta > 1$, either $\delta = \delta_5 = 1$ or $\delta = \delta_8 = 1$ and
$\underset{k=1}{\overset{N}\sum} \; |\alpha_k| \le \; 1 $.\\
5. $\Theta > 1$, $\delta = \delta_4 = 1$ and 
$\underset{k=1}{\overset{N}\sum} \;|\gamma_k| \le \; 1$.
\end{corl}
\begin{Rem}
1. In the critical case $\Gamma = 1,$ the fractal dimension bounds
of  CHFIF $f_1(x)$  found in (\ref{fd1}) coincide with
the fractal dimension bounds of FIF $f_2(x)$ found  in \cite{G}, if $f_1(x)$ is also self-affine.\\
2. Choosing  $\alpha_i$ suitably, $ \Omega = \Theta = 1$
for any self-affine CHFIF $f_1 = f_2$. In the resulting critical cases, 
the bounds on fractal dimension of any self-affine FIF $f_1$ can be found by (\ref{fd})
with suitable choice of hidden variables even if $\Gamma \ne 1$. \\
3. In Corollary \ref{cor2}, critical  CHFIF $f_1(x)$ is
considered as fractal function, since $f_1(x)$ satisfies
$\omega (f_1,x) = \bigcirc ( |x|^{\delta} \log|x|).$
Consequently, fractal functions having $ D_B(graph(f_1)) = 1$ can be constructed
by using this corollary.
\end{Rem}
\section{EXAMPLES}
Consider the interpolation data \{(0,2),(0.35,7),(.75,4),(1,9)\}. Here,
for simplicity, we construct affine  CHFIFs.
Since, in this case $ \lambda = \mu = 1,$  it follows that $\Theta = \Omega.$
In Figs. 1-3, the generalized set of data is chosen such that $z_i = y_i$  and
in Figs. 4-16, the generalized set of data chosen such that
$z_i \ne y_i.$ The  values of $ \alpha_i, \beta_i \;  \text{and} \; \gamma_i $
chosen  for the computer generation of affine  CHFIFs for all these
figures are given in Table 1. Fig. 1 gives the self-affine  CHFIF $f_1(x)$
for the given interpolation data whenever $ \alpha_i +  \beta_i =  \gamma_i $.
Fig. 2 and Fig. 3 show respectively the effect on the CHFIF for suitable choices of
$ \alpha_i,  \beta_i$ and  $ \gamma_i$  when the effective scaling factor is 
close to $-2^+$ and $2^-$.
\begin{center}
{{\bf Table 1 :}{ \it Free variables and constrained free variables in the
construction of affine CHFIFs}}
\end{center}
\begin{center}
\begin{tabular}{|c|c|c|c|c|c|c|c|c|c|}
\hline {\bf Figures} &  $\mathbf {\alpha_1}$ & $\mathbf {\alpha_2}$ & $\mathbf {\alpha_3}$
 & $\mathbf {\beta_1}$ & $\mathbf {\beta_2}$ & $ \mathbf {\beta_3}$ &$ \mathbf {\gamma_1}$ &
 $\mathbf {\gamma_2}$  & $\mathbf {\gamma_3}$ \\ \hline
 1 & 0.8 & 0.7 & 0.3 & -0.3 & -0.4 & -0.2 & 0.5 & 0.3 & 0.6 \\ \hline
 2 & 0.99 & 0.99 & 0.99 & 0.99 & 0.99 & 0.99 & 0.005 & 0.005 & 0.005 \\ \hline
 3 & -0.999 &-0.999 & -0.999 & -0.99 & -0.99 & -0.005 & -0.005 & -0.005 & -0.005 \\ \hline
 4 & 0.2 & 0.38 & 0.2 & 0.4 & 0.35 & 0.5 &  0.3 & 0.3 & 0.24 \\ \hline
 5 & 0.2 & 0.4 & 0.22 & 0.4 & 0.35 & 0.5 &  0.35 & 0.3 & 0.2 \\ \hline
 6 & 0.4 & 0.3 & 0.5 & 0.4 & 0.35 & 0.5 &  0.3 & 0.5 & 0.4 \\ \hline
 7 & 0.2 & 0.38 & 0.2 & 0.4 & 0.35 & 0.5 &  0.3 & 0.5 & 0.4 \\ \hline
 8 & 0.4 & 0.3 & 0.5 & 0.4 & 0.35 & 0.5 &  0.3 & 0.3 & 0.24 \\ \hline
 9 & 0.2 & 0.4 & 0.22 & 0.4 & 0.35 & 0.5 &  0.3 & 0.3 & 0.24 \\ \hline
 10 & 0.2 & 0.38 & 0.2 & 0.4 & 0.35 & 0.5 &  0.35 & 0.3 & 0.2 \\ \hline
11 & 0.2 & 0.4 & 0.22 & 0.4 & 0.35 & 0.5 &  0.3 & 0.5 & 0.4 \\ \hline
12 & 0.4 & 0.3 & 0.5 & 0.4 & 0.35 & 0.5 &  0.35 & 0.3 & 0.2 \\ \hline
 13 & 0.2 & 0.38 & 0.2 & -0.6 & -0.45 & -0.4 &  0.3 & 0.3 & 0.24 \\ \hline
 14 & 0.4 & 0.3 & 0.5 & -0.6 & -0.45 & -0.4 &  0.3 & 0.5 & 0.4 \\ \hline
 15 & 0.4 & 0.3 & 0.5 & -0.6 & -0.45 & -0.4 &  0.3 & 0.3 & 0.24 \\ \hline
 16 & 0.4 & 0.3 & 0.5 & -0.6 & -0.45 & -0.4 &  0.3 & 0.5 & 0.4 \\ \hline
\end{tabular}
\end{center}
\par
Figs. 4-12  with the  fixed values of $\beta_i$ ( as in Table 1 ) and $z_i$ 
( 3,1,8,5 respectively for interpolation data points ) illustrate the
nature of non-self-affine  CHFIF, depending upon
various cases of smoothness analysis in Theorems \ref{th1}-\ref{th3}, when $\Theta = \Omega.$
As expected,  CHFIFs in these figures have the same type 
of shape  since the underlying function spaces are independent of
$\beta_i$ and $z_i.$  Figs. 13-15 give the effect of change in the
values of $\beta_i$ ( as in Table 1) on the shape of non-self-affine
CHFIF. Fig. 16 shows the effect of change in
the value of $z_i$ ( 7,9,10,8 respectively for interpolation data points )
on the shape of non-self-affine  CHFIF.
On comparing Fig. 4 with Fig. 13, Fig. 6 with Fig. 14 and  Fig. 11 with
Fig. 15, it is found that although  CHFIFs
are in the same function spaces, these are very much different
in shape due to  changes in the values of $\beta_i$ ( as given
in Table 1 ). The underlying function spaces are the same
because these spaces depend only  on the values of $\alpha_i, \gamma_i,
\lambda_i$  and $\mu_i.$ Further, comparing Fig. 14 with Fig. 16, it is
observed that by keeping all the other values  fixed and changing 
only the values of $z_i$ from 3,1,8,5 respectively to 7,9,10,8 
in generalized  interpolation data, the shape of the CHFIF changes arbitrarily.
\section{CONCLUSION}
A generalized IFS is constructed in the present paper 
for generating coalescence  hidden variable FIF.  
The existence and uniqueness of the CHFIF
is proved by choosing suitable values of the variables $\alpha_i, \beta_i$ and $\gamma_i$
and the parameter $z_i$. Our IFS gives  CHFIFs that may be self-affine or
non-self-affine depending on free variables, constraints free variable and the parameters $z_i$..
When construction of  the CHFIF is carried out  by adding $n$ dimensions
linearly in generalized interpolation data, 
{\it{ $(n+1)$ free variables}} and at most {\it {$(1+2+ \dots +n)$ constrained
free variables}} can be chosen. If all of the extra $n$ dimensions take the same values of $ z_i,$
the scaling factor of  the CHFIF lies between $-(n+1)^+$ and $(n+1)^-.$
Besides using the generalized  IFS for construction of CHFIFs in the present work,
it can also be used in other scientific
applications to capture the self-affine and non-self-affine  nature simultaneously for the
relevant curves.
\par
It is seen that the  smoothness of  CHFIF $f_1(x)$
depends on free variables $\alpha_i$ and $\gamma_i$ as well as on the smoothness of $p_i(x)$ and
$q_i(x).$ Although, $z_i$ and $\beta_i$ are responsible for the
shape of  the CHFIF, these  are found  not to  affect its smoothness.
In general, the deterministic construction of functions having order of
modulus of continuity  $O(|t|^{\delta}(\log|t|)^m)$ ( $m$ a non-negative integer,
and $ 0 < \delta \le 1$) is possible through  the CHFIF.
The fact that CHFIFs are different in shape although they are in the same
function spaces may enable considering them in more general function
spaces such as Besov and Triebel-Lizorkin spaces apart from Lipschitz
spaces. These former spaces have additional indices that `fine-tune'
a function. Our bounds of fractal dimension of  CHFIFs are found
in different critical conditions.  Finally, it is proved that by
suitable choices of the hidden variables, the fractal dimension bounds 
for any self-affine FIF can be found using the bounds obtained with
 the critical condition $\Omega = \Theta = 1 $.

\newpage
\clearpage
\begin{figure}
\begin{minipage}{0.45\textwidth}
\psfig{file=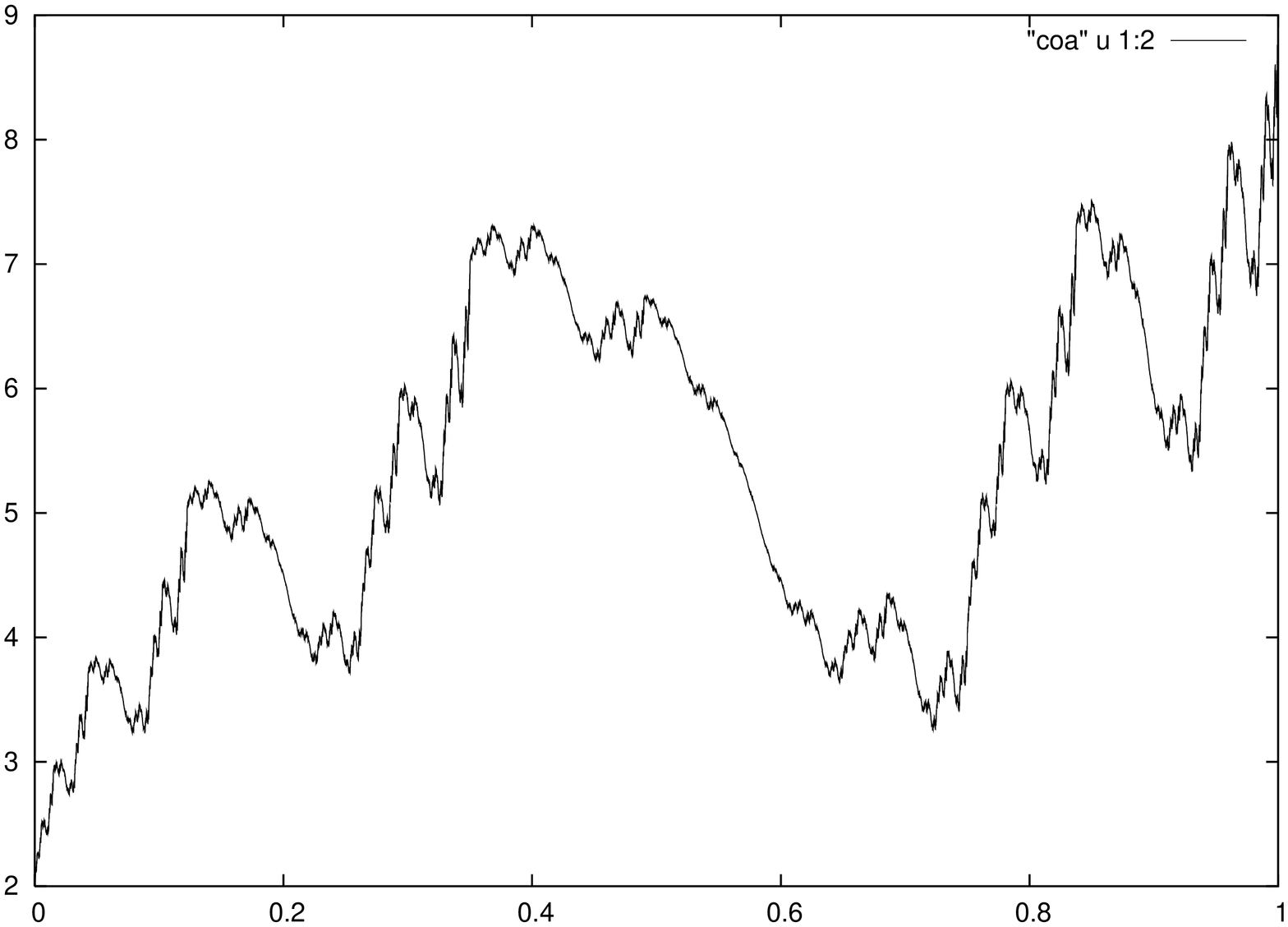,scale=0.25,angle=-90} \\
\centering{\small Fig. 1 Self-affine  CHFIF $f_1(x)$.}
\end{minipage} \hfill
\begin{minipage}{0.45\textwidth}
\psfig{file=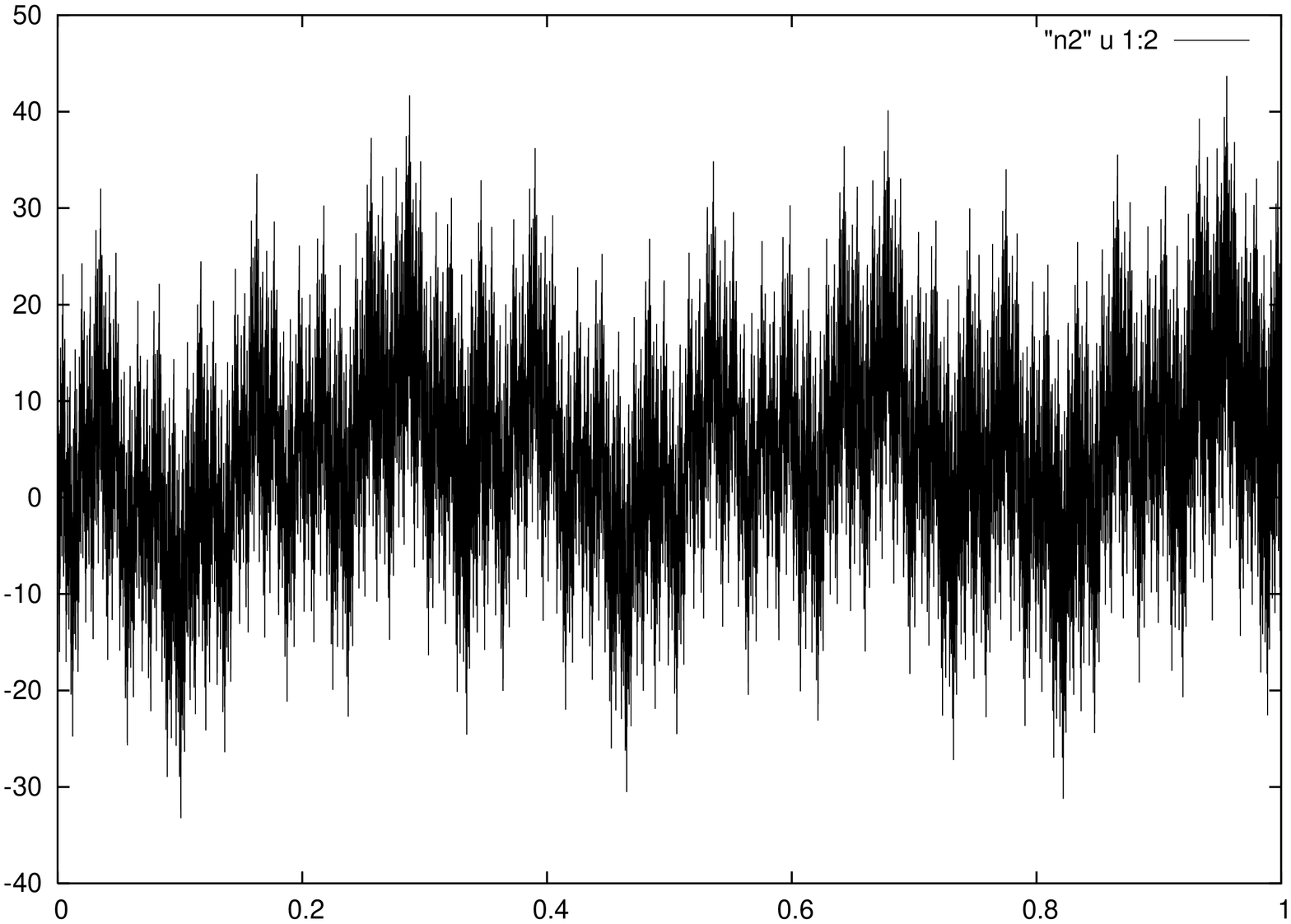,scale=0.25,angle=-90}\\
\centering{\small Fig. 2  CHFIF $f_1(x)$ with  scaling factor $-2^+$}
\end{minipage}\\
\begin{minipage}{0.45\textwidth}
\psfig{file=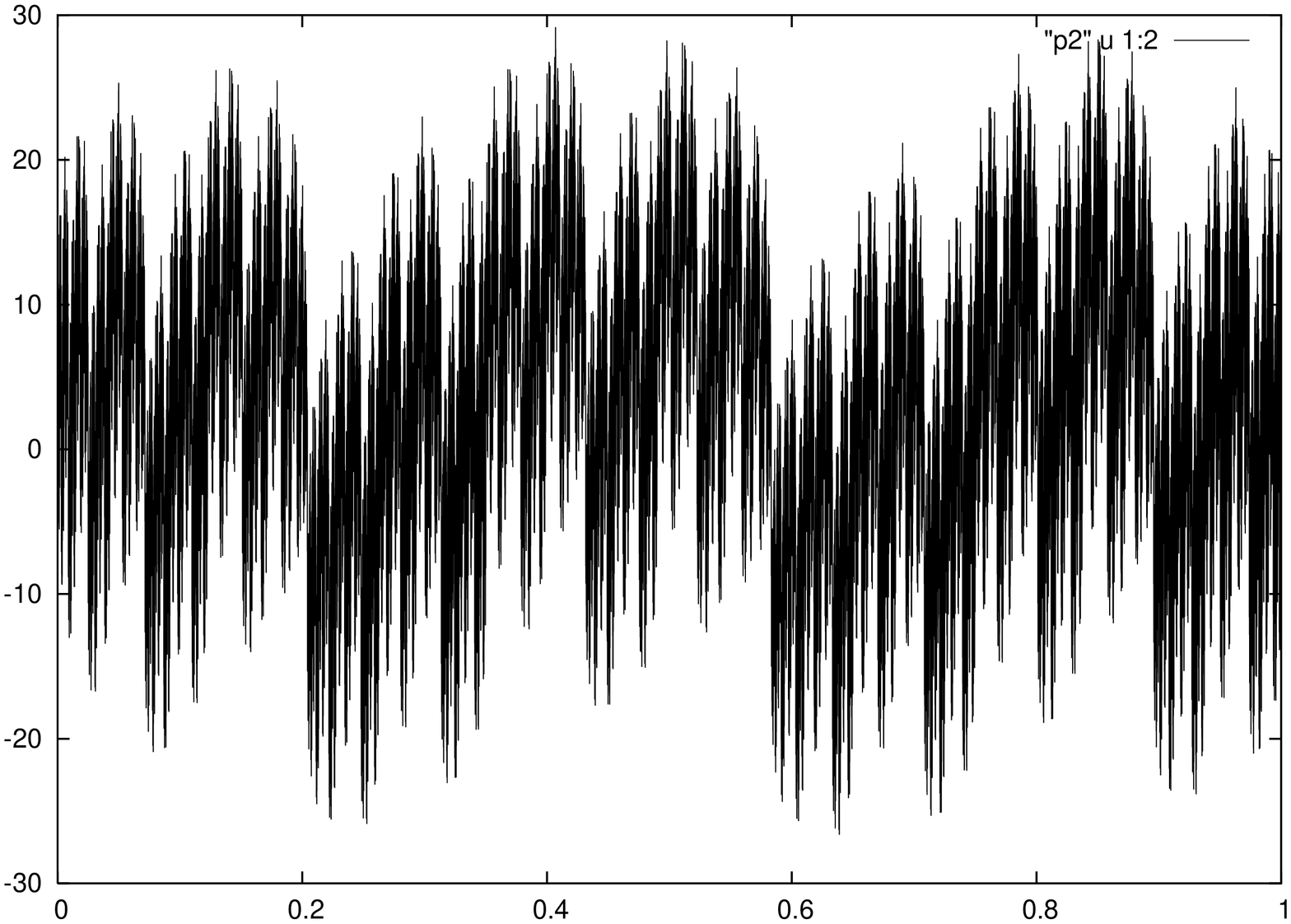,scale=0.25,angle=-90}\\
\centering{\small Fig. 3   CHFIF  $f_1(x)$ with  scaling factor $2^-$}
\end{minipage} \hfill
\begin{minipage}{0.45\textwidth}
\psfig{file=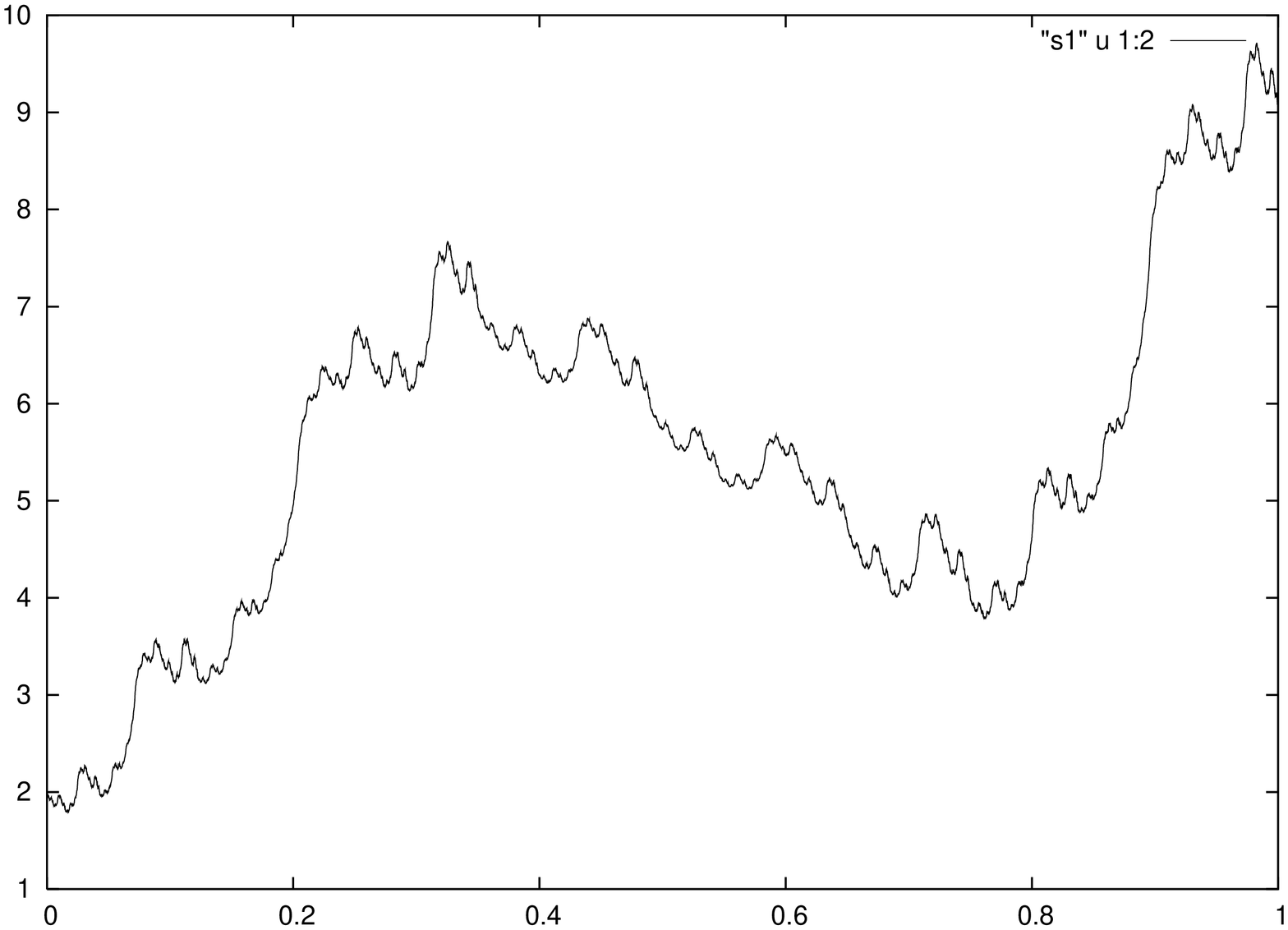,scale=0.25,angle=-90} \\
\centering{\small Fig. 4  CHFIF with $\Theta = \Omega < 1 , \Gamma < 1 $.}
\end{minipage}\\
\begin{minipage}{0.45\textwidth}
\psfig{file=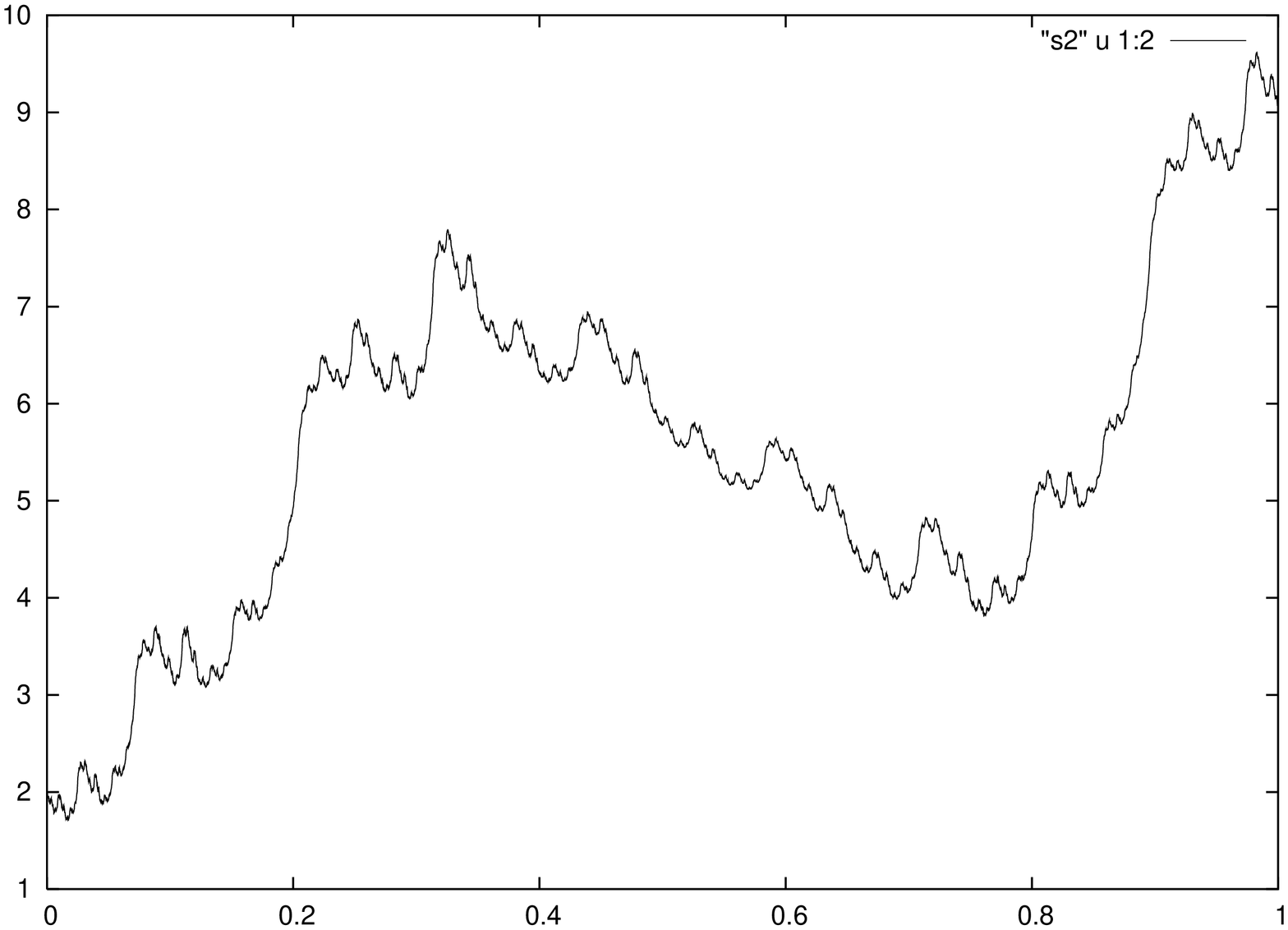,scale=0.25,angle=-90}\\
\centering{\small Fig. 5  CHFIF with $\Theta = \Omega = 1 , \Gamma = 1 $}
\end{minipage} \hfill
\begin{minipage}{0.45\textwidth}
\psfig{file=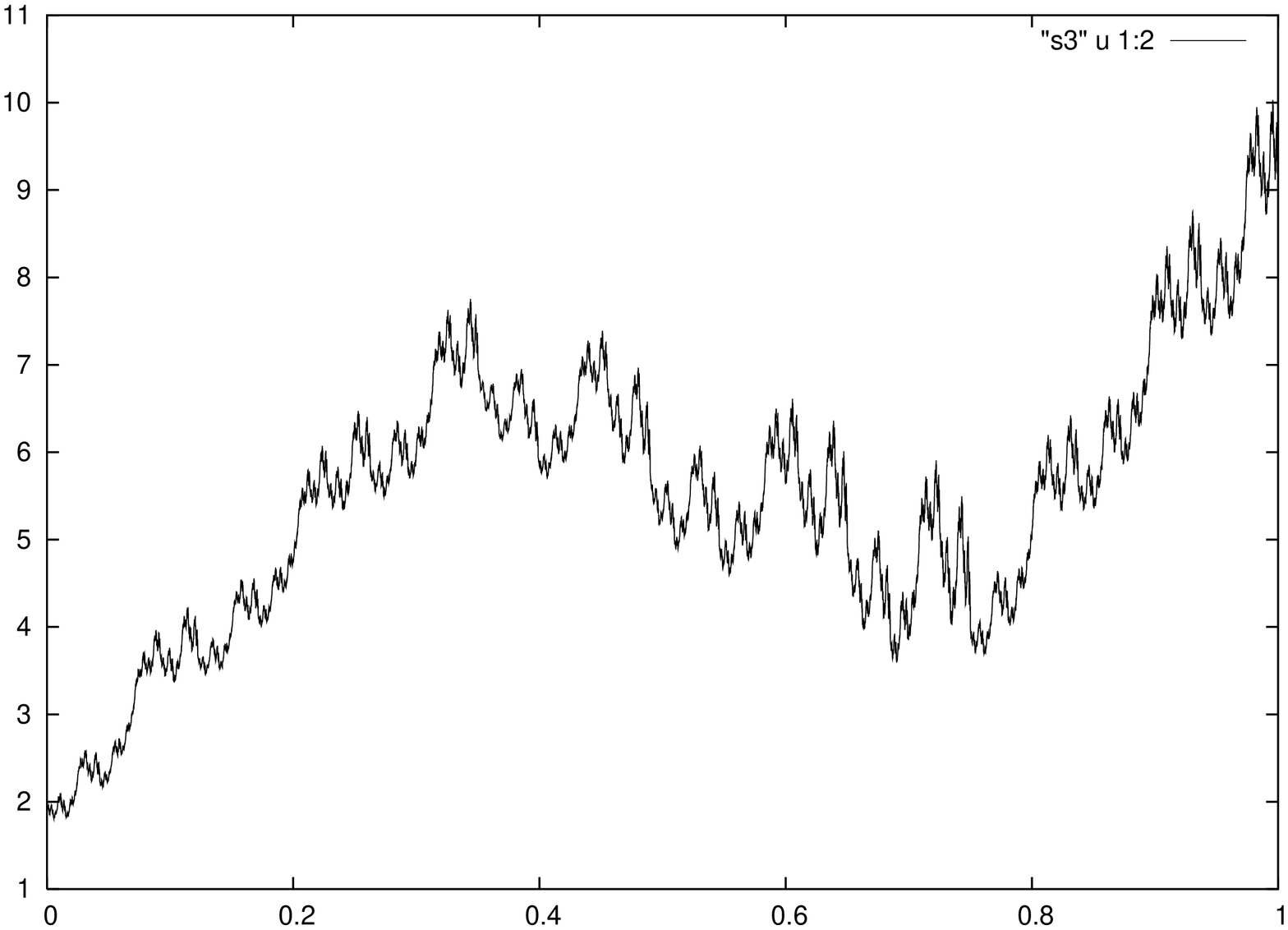,scale=0.25,angle=-90}\\
\centering{\small Fig. 6  CHFIF with $\Theta = \Omega > 1 , \Gamma > 1 $}
\end{minipage} \\
\begin{minipage}{0.45\textwidth}
\psfig{file=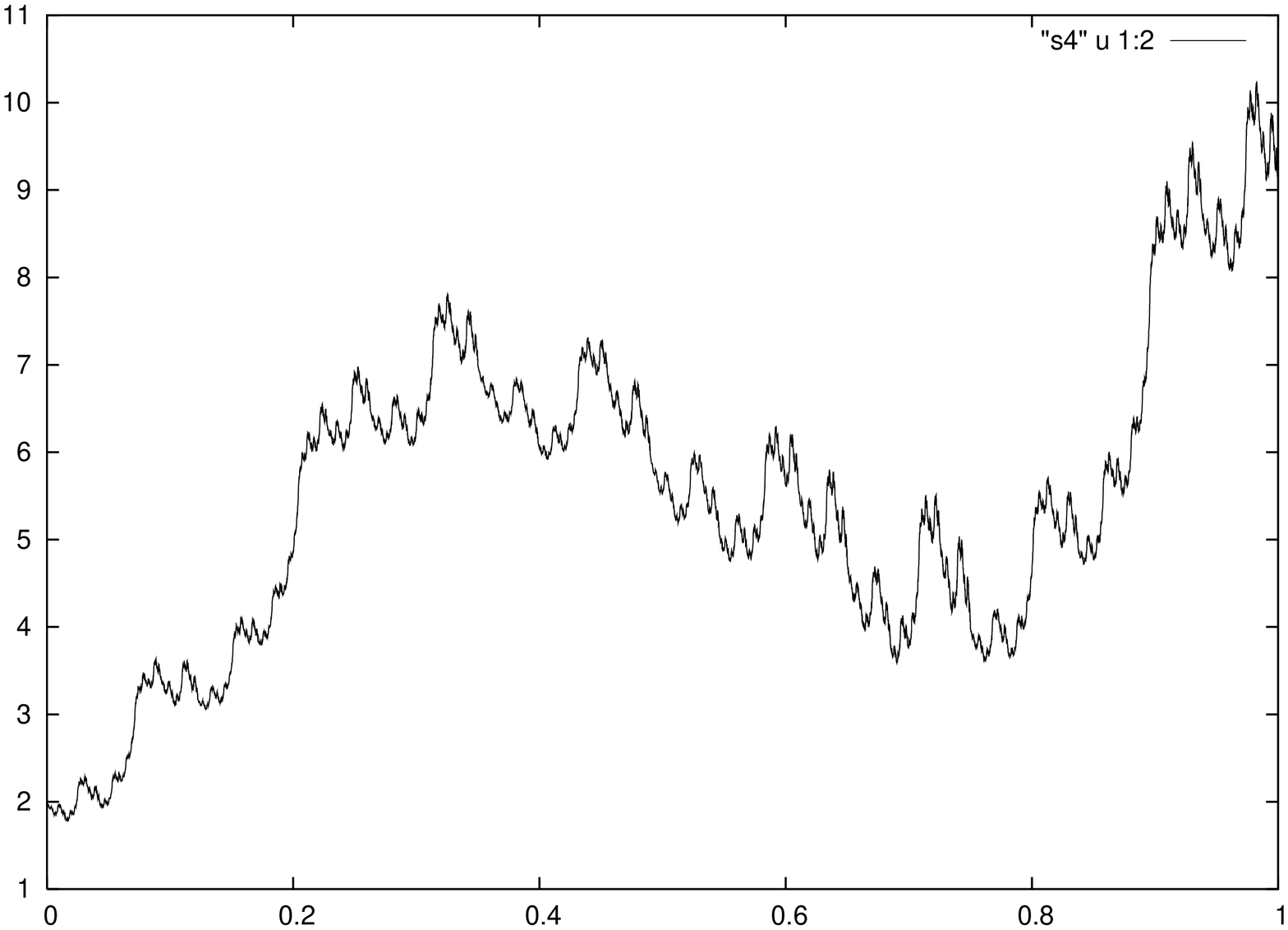,scale=0.25,angle=-90}\\
\centering{\small Fig. 7  CHFIF with $ \Theta = \Omega < 1 , \Gamma > 1 $}
\end{minipage}\hfill
\begin{minipage}{0.45\textwidth}
\psfig{file=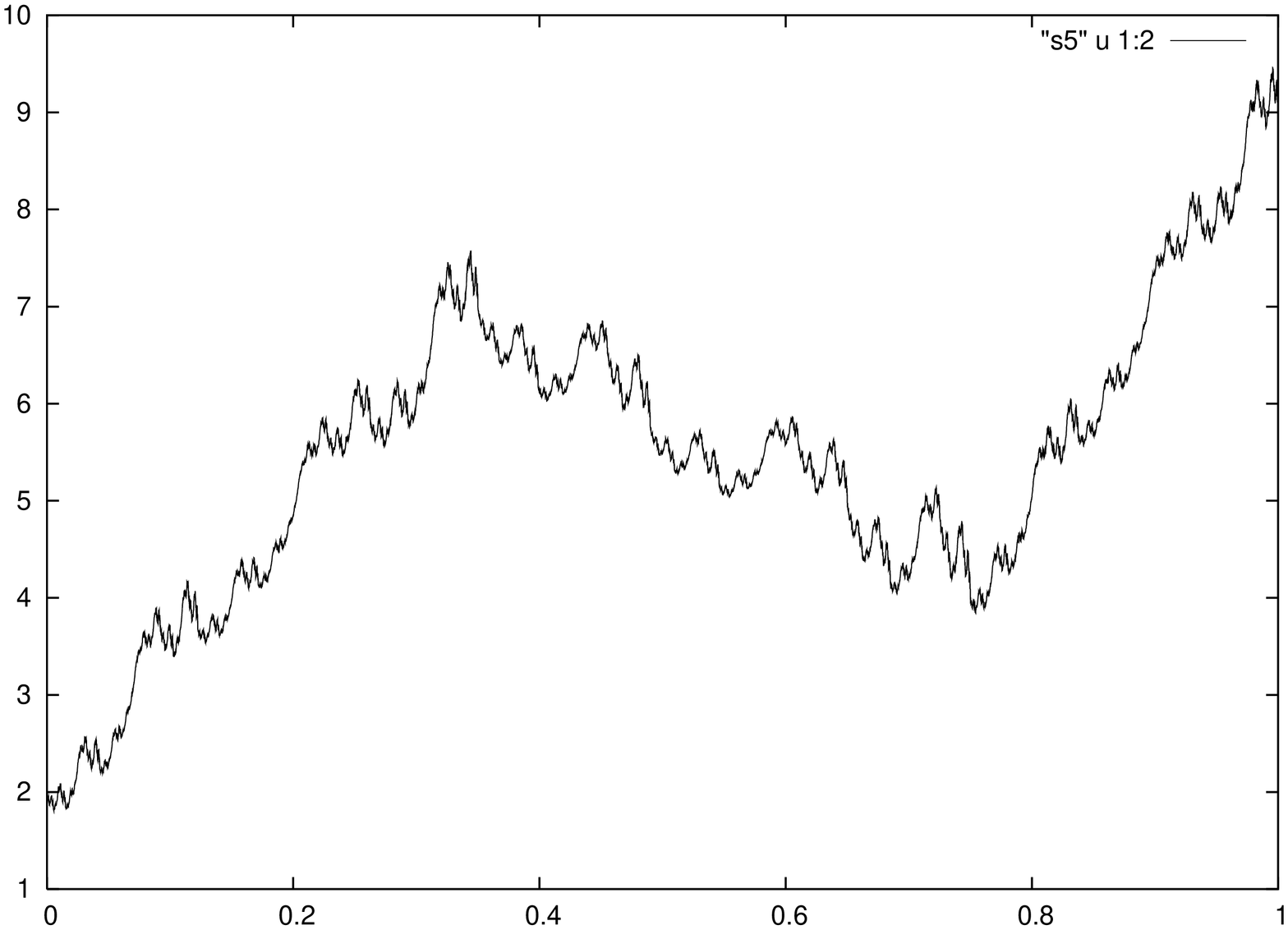,scale=0.25,angle=-90} \\
\centering{\small Fig. 8  CHFIF with $\Theta = \Omega > 1 , \Gamma < 1 $.}
\end{minipage}
\end{figure}

\newpage
\clearpage
\begin{figure}
\begin{minipage}{0.45\textwidth}
\psfig{file=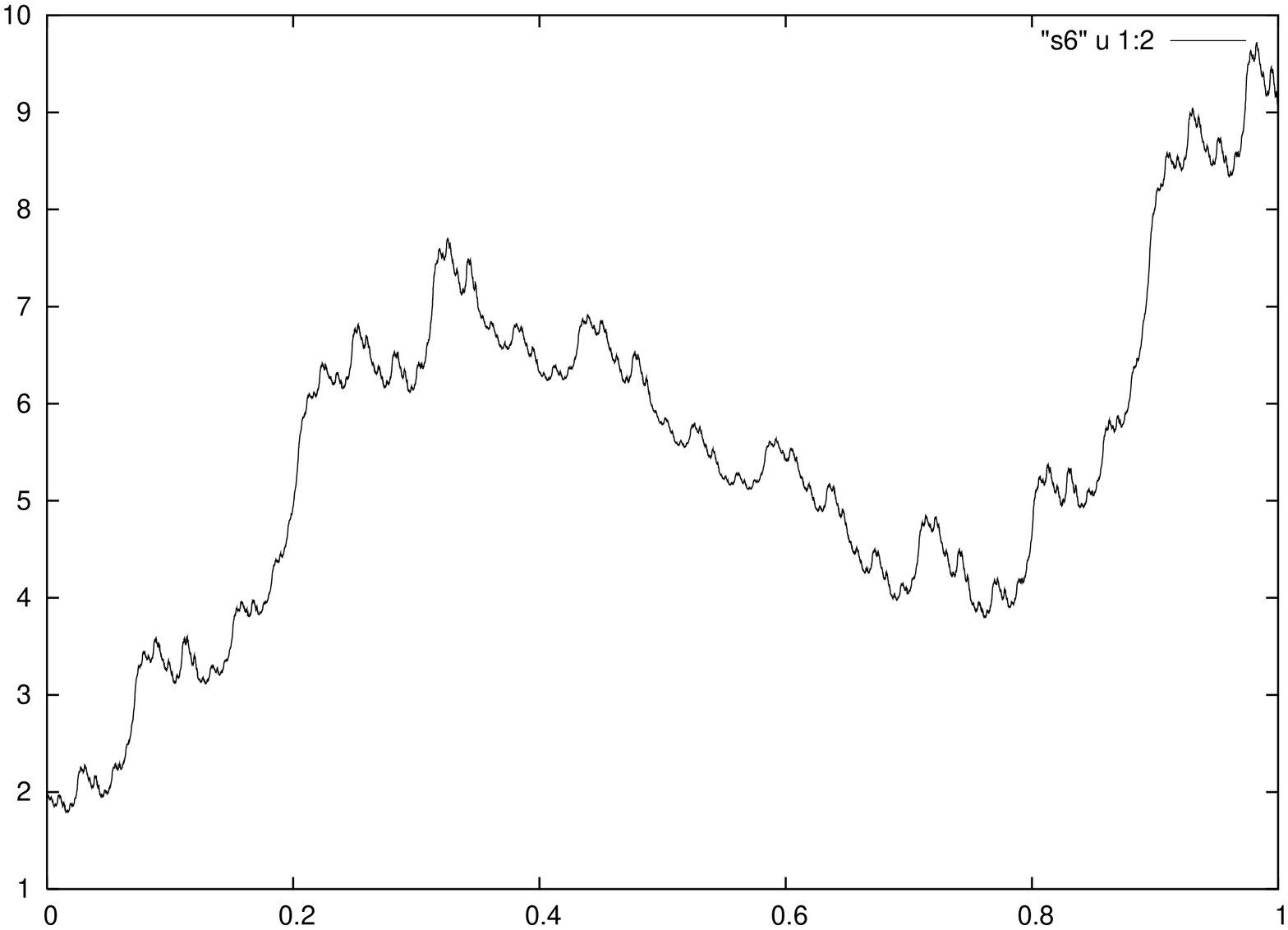,scale=0.25,angle=-90}\\
\centering{\small Fig. 9  CHFIF with $ \Theta = \Omega = 1 , \Gamma < 1 $}
\end{minipage} \hfill
\begin{minipage}{0.45\textwidth}
\psfig{file=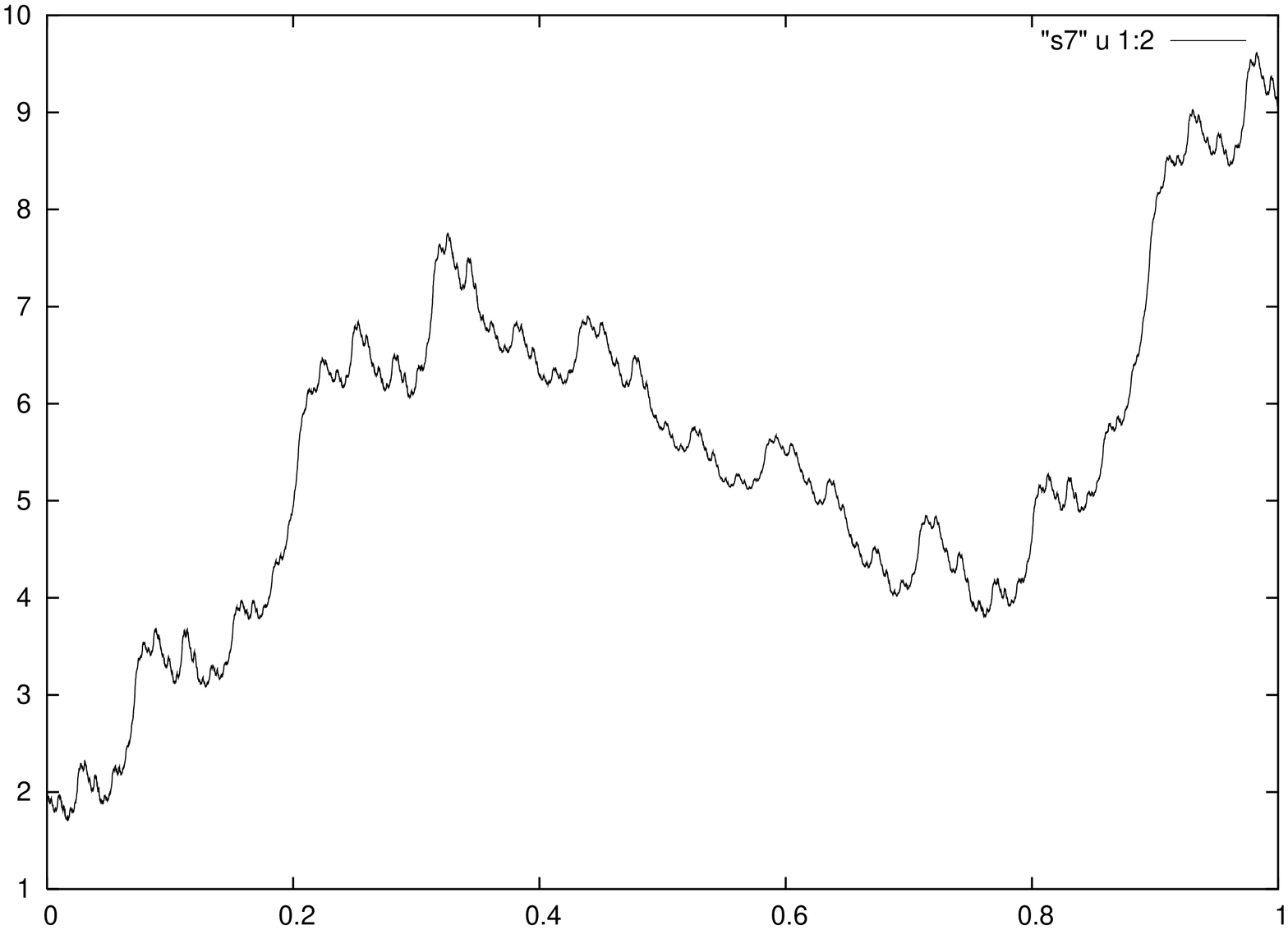,scale=0.25,angle=-90}\\
\centering{\small Fig. 10  CHFIF with $ \Theta = \Omega < 1 , \Gamma = 1 $}
\end{minipage}\\
\begin{minipage}{0.45\textwidth}
\psfig{file=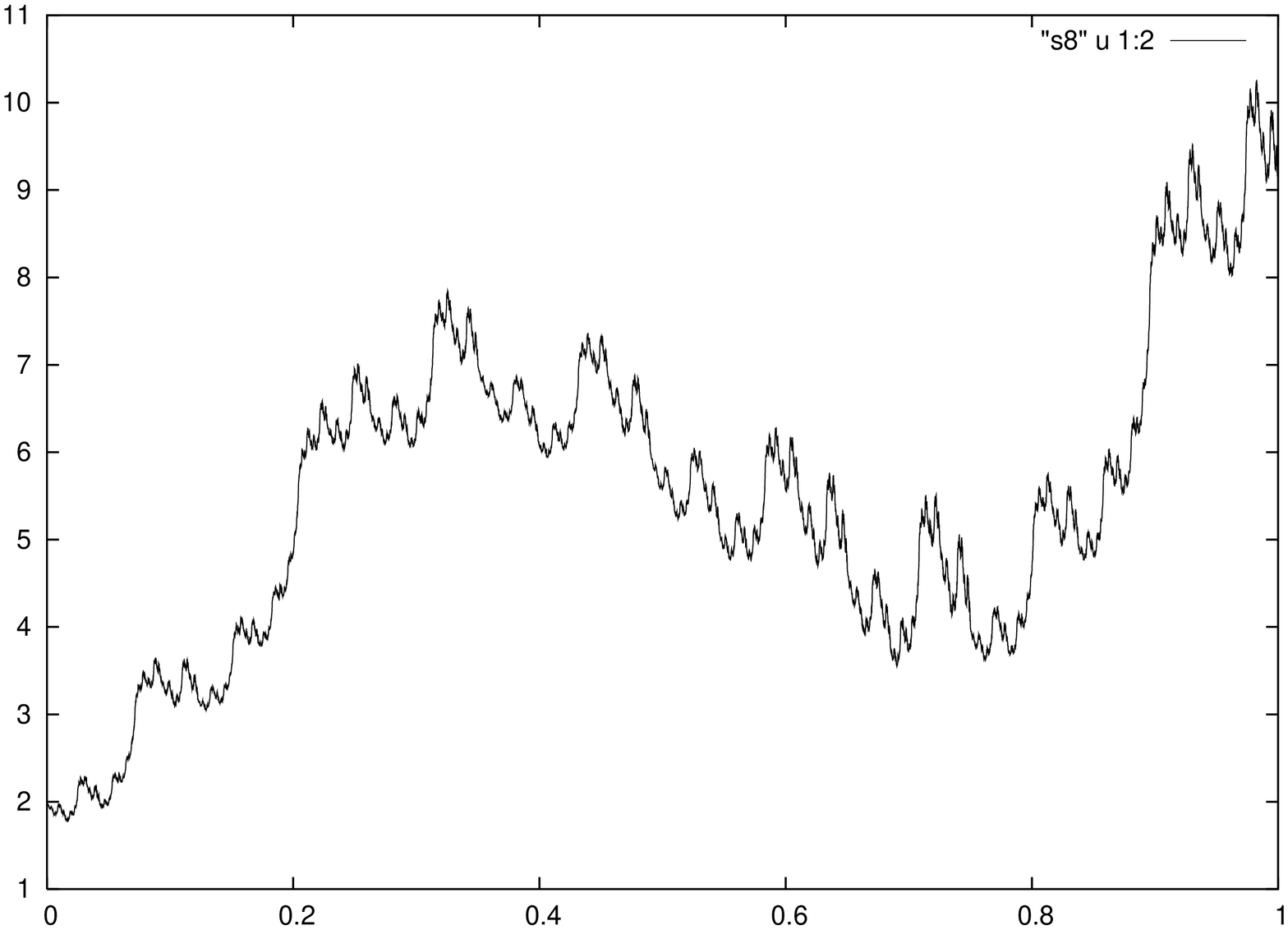,scale=0.25,angle=-90}\\
\centering{\small Fig. 11  CHFIF with $ \Theta = \Omega = 1 , \Gamma > 1 $}
\end{minipage} \hfill
\begin{minipage}{0.45\textwidth}
\psfig{file=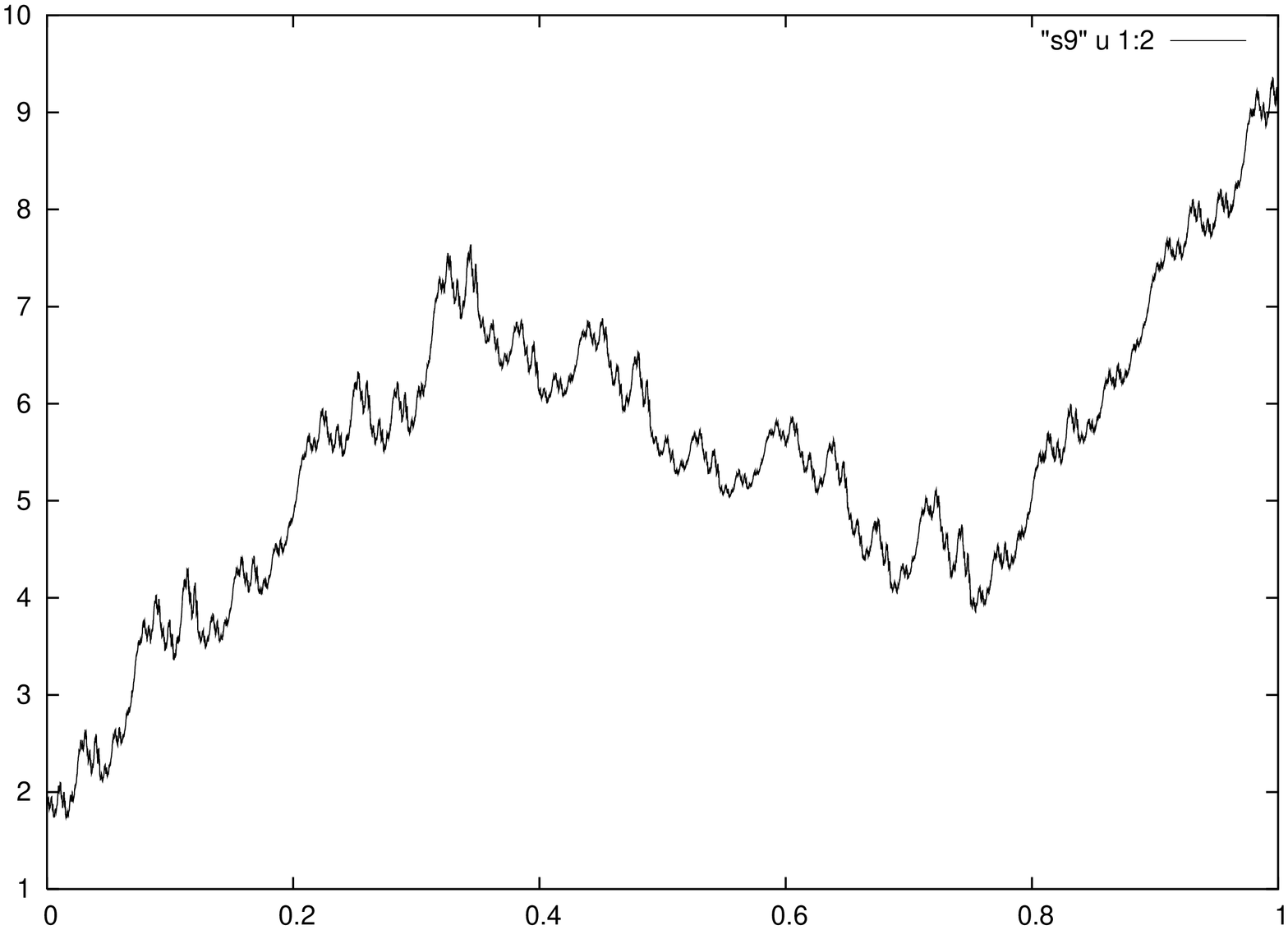,scale=0.25,angle=-90}\\
\centering{\small Fig. 12  CHFIF with $ \Theta = \Omega > 1 , \Gamma = 1 $}
\end{minipage}\\
\begin{minipage}{0.45\textwidth}
\psfig{file=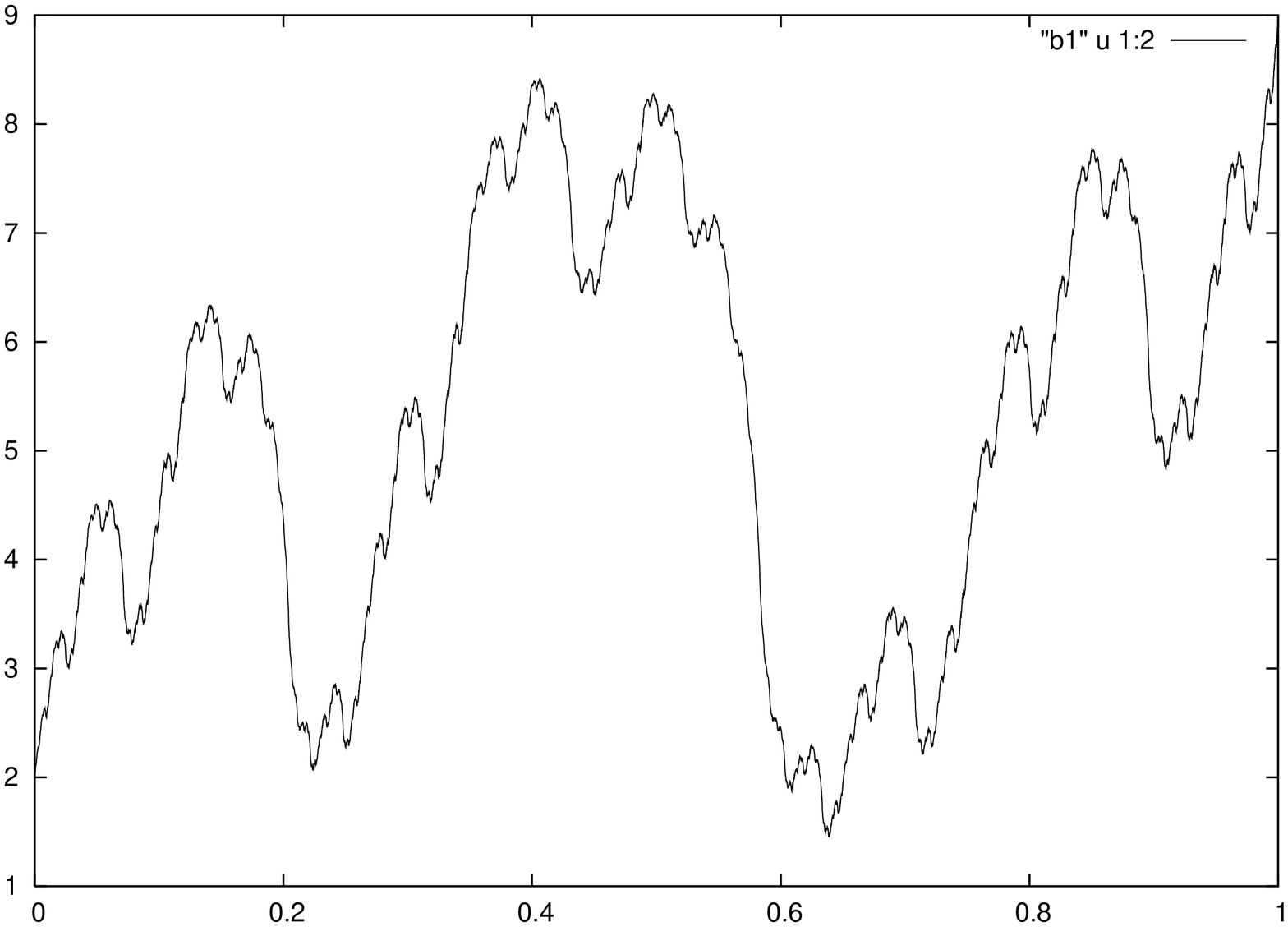,scale=0.25,angle=-90}\\
\centering{\small Fig. 13  CHFIF with $ \Theta = \Omega < 1 , \Gamma < 1 $
with a different set of $\beta_i.$}
\end{minipage} \hfill
\begin{minipage}{0.45\textwidth}
\psfig{file=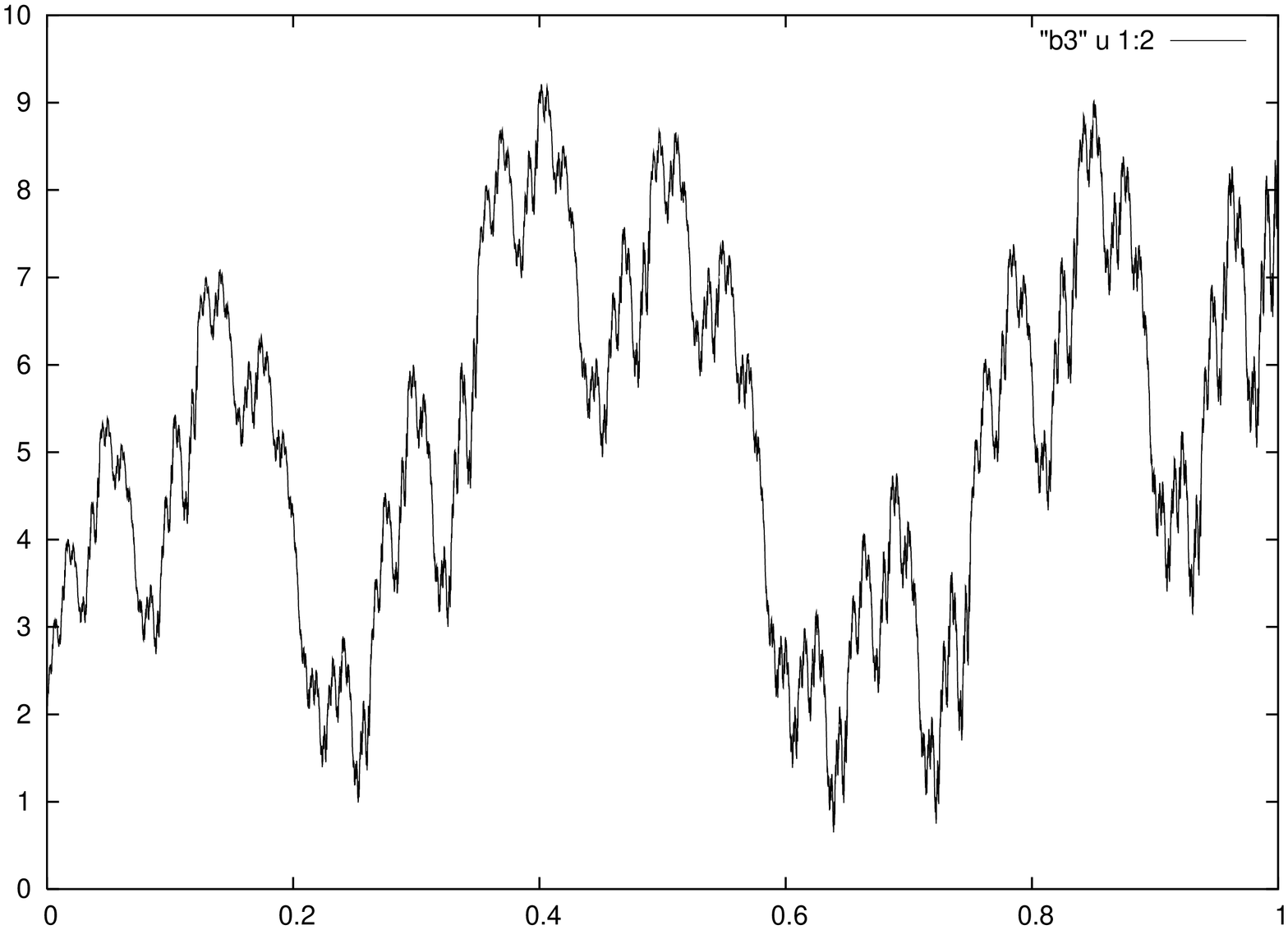,scale=0.25,angle=-90}\\
\centering{\small Fig. 14  CHFIF with $ \Theta = \Omega > 1 , \Gamma > 1 $
with a different set of $\beta_i.$}
\end{minipage}\\
\begin{minipage}{0.45\textwidth}
\psfig{file=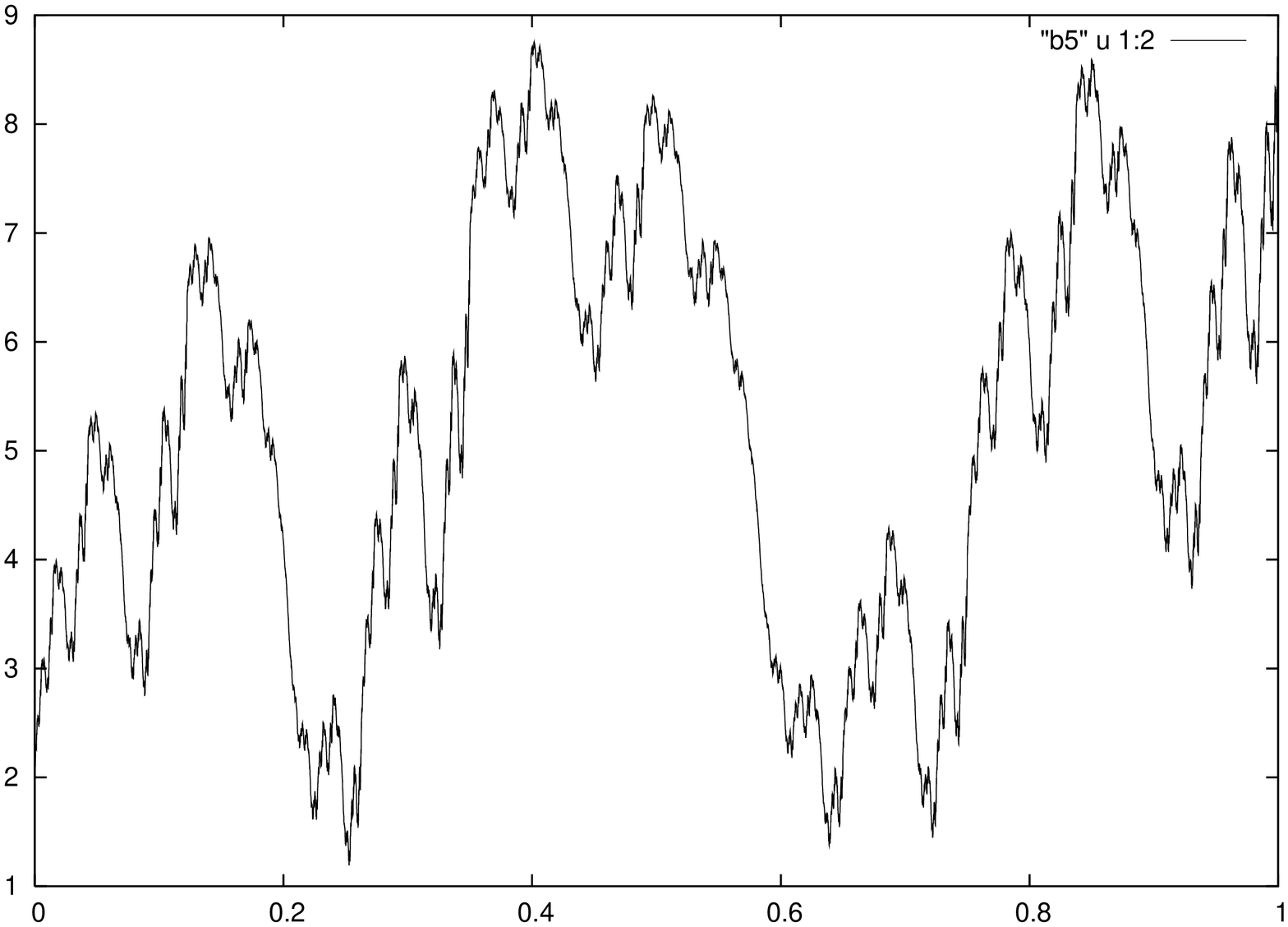,scale=0.25,angle=-90}\\
\centering{\small Fig. 15   CHFIF with $ \Theta = \Omega > 1 , \Gamma < 1 $
with a different set of $\beta_i.$}
\end{minipage} \hfill
\begin{minipage}{0.45\textwidth}
\psfig{file=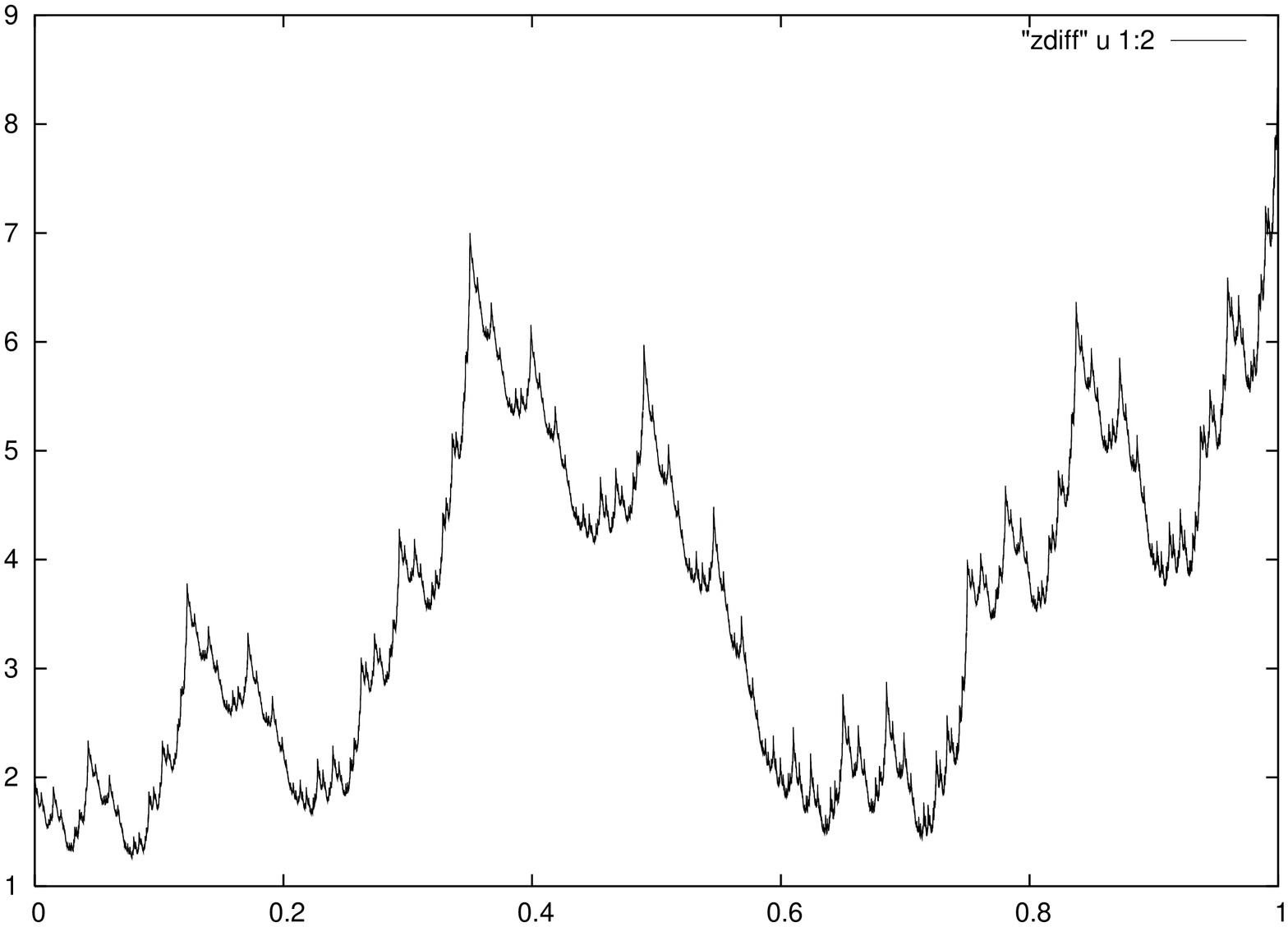,scale=0.25,angle=-90}\\
\centering{\small Fig. 16  CHFIF with $ \Theta = \Omega > 1 , \Gamma > 1 $
with a different set of $z_i.$}
\end{minipage}
\end{figure}
\end{document}